\documentclass[journal]{IEEEtran}
\pdfoutput=1

\usepackage{ifpdf}
\usepackage{xifthen}
\usepackage[l2tabu, orthodox]{nag}

\ifCLASSINFOpdf
  \usepackage[pdftex]{graphicx}
   \DeclareGraphicsExtensions{.pdf,.jpg,.jpeg,.png}
\else
  \usepackage[dvips]{graphicx}
   \DeclareGraphicsExtensions{.eps}
\fi

\usepackage[cmex10]{amsmath}
\usepackage{mathtools,amssymb,amsthm}
\usepackage[norelsize,vlined,ruled]{algorithm2e}		% Only use with option [H]!
\usepackage{xfrac}				% Needed for command \sfrac

\usepackage{pgfplots}
\pgfplotsset{compat=1.9}
\usepackage{xcolor}
\usepackage{tikz}
\usetikzlibrary{calc,decorations.shapes,shapes.geometric}
\usepgfplotslibrary{statistics}
\usepackage{dblfloatfix}

\usepackage{booktabs}			% Improved tables
\usepackage{array}
\newcolumntype{h}{c<{\hskip0.43cm}}%
\newcolumntype{H}{c<{\hskip0.63cm}}%

\usepackage{fix-cm}				% Fix error with missing font sizes
\usepackage{eurosym}
\usepackage[final]{microtype}	% Improved typesetting

\usepackage{url}
\usepackage{cite}
\usepackage[absolute]{textpos}
\DontPrintSemicolon

\makeatletter
\newcommand{\removelatexerror}{\let\@latex@error\@gobble}
\makeatother

\renewcommand{\exp}[1]{\text{e}^{#1}}

\newcommand{\set}[1]{\left\{#1\right\}}

\newcommand{\bigO}{{\mathcal{O}}}

\newcommand{\squeeze}[2][0.25]{
	{\text{% Interrupt math to shrink spaces
		\thinmuskip=#1\thinmuskip%
		\medmuskip=#1\medmuskip%
		\thickmuskip=#1\thickmuskip%
		$#2$%
	}}%
}

\newcommand{\eg}{e.g.\ }
\newcommand{\cf}{c.f.\ }
\newcommand{\ie}{i.e.\ }

\newcommand{\enquote}[1]{``#1''}

\makeatletter
\def\discrange{\@ifnextchar*\discrangeStar{\discrangeInternal{n}}}
\makeatother
\def\discrangeStar*{\discrangeInternal{y}}
\def\discrangeInternal#1#2#3{{%
	\lbrack%
	\ifthenelse{\equal{#1}{y}}{%
		\squeeze[0.1]{#2}\:..\:\squeeze[0.1]{#3}
	}{%
		#2\:..\:#3
	}%
	\rbrack%
}}

\makeatletter
\newcommand{\renewvar}[3]{%
	\expandafter\def\csname #1\endcsname{\@ifnextchar({\csname #1Parentheses\endcsname}{\csname #1Braces\endcsname}}%
	\ifcase #2
		\expandafter\def\csname #1\endcsname{#3\@ifnextchar({\PackageError{newvar}{The variable #1 should have no parameters!}{}}\relax}%
		\expandafter\def\csname #1Braces\endcsname{#3}%
	\or	\expandafter\def\csname #1Parentheses\endcsname(##1){#3}%
		\expandafter\def\csname #1Braces\endcsname##1{#3}%
	\or	\expandafter\def\csname #1Parentheses\endcsname(##1,##2){#3}%
		\expandafter\def\csname #1Braces\endcsname##1##2{#3}%
	\or	\expandafter\def\csname #1Parentheses\endcsname(##1,##2,##3){#3}%
		\expandafter\def\csname #1Braces\endcsname##1##2##3{#3}%
	\or	\expandafter\def\csname #1Parentheses\endcsname(##1,##2,##3,##4){#3}%
		\expandafter\def\csname #1Braces\endcsname##1##2##3##4{#3}%
	\or	\expandafter\def\csname #1Parentheses\endcsname(##1,##2,##3,##4,##5){#3}%
		\expandafter\def\csname #1Braces\endcsname##1##2##3##4##5{#3}%
	\or	\expandafter\def\csname #1Parentheses\endcsname(##1,##2,##3,##4,##5,##6){#3}%
		\expandafter\def\csname #1Braces\endcsname##1##2##3##4##5##6{#3}%
	\or	\expandafter\def\csname #1Parentheses\endcsname(##1,##2,##3,##4,##5,##6,##7){#3}%
		\expandafter\def\csname #1Braces\endcsname##1##2##3##4##5##6##7{#3}%
	\or	\expandafter\def\csname #1Parentheses\endcsname(##1,##2,##3,##4,##5,##6,##7,##8){#3}%
		\expandafter\def\csname #1Braces\endcsname##1##2##3##4##5##6##7##8{#3}%
	\or	\expandafter\def\csname #1Parentheses\endcsname(##1,##2,##3,##4,##5,##6,##7,##8,##9){#3}%
		\expandafter\def\csname #1Braces\endcsname##1##2##3##4##5##6##7##8##9{#3}%
	\fi
}
\newcommand{\newvar}[3]{
	\@ifundefined{#1}{}{\PackageError{newvar}{The variable #1 is already exists!}{Use renewvar to redefine existing variables.}}
	\renewvar{#1}{#2}{#3}
}
\makeatother

\newvar{leftSizeT}{2}{\leftSizeSymbol^{#1}\ifthenelse{\isempty{#2}}{}{(#2)}}
\newvar{leftSize}{1}{\leftSizeSymbol(#1)}
\newcommand{\rightSizeSymbol}{r}
\newvar{rightSizeT}{2}{\rightSizeSymbol^{#1}\ifthenelse{\isempty{#2}}{}{(#2)}}
\newvar{rightSize}{1}{\rightSizeSymbol(#1)}
\newcommand{\sizeSymbol}{s}
\newvar{sizeT}{2}{\sizeSymbol_{#1}\ifthenelse{\isempty{#2}}{}{(#2)}}
\newvar{size}{1}{\sizeSymbol(#1)}
\newvar{leftChildT}{2}{\text{llink}^{#1}\ifthenelse{\isempty{#2}}{}{\left({#2}\right)}}
\newvar{rightChildT}{2}{\text{rlink}^{#1}\ifthenelse{\isempty{#2}}{}{\left({#2}\right)}}
\newvar{leftChild}{1}{\leftChildT{}{#1}}
\newvar{rightChild}{1}{\rightChildT{}{#1}}
\newvar{topleftT}{2}{\mu^{#1}_{#2}}
\newvar{topleft}{1}{\topleftT{}{#1}}
\newvar{toprightT}{2}{\eta^{#1}_{#2}}
\newvar{topright}{1}{\toprightT{}{#1}}

\definecolor{TUMblue}{RGB}{0,101,189} % #0065bd, HSB 208,1,0.37
\definecolor{TUMivory}{RGB}{218,215,203}
\definecolor{TUMorange}{RGB}{227,114,34}

\definecolor{TUMred}{RGB}{189,0,6}
\definecolor{TUMyellow}{RGB}{189,183,0}

\definecolor{TUMgreen}{RGB}{101,189,0}
\definecolor{TUMmagenta}{RGB}{189,0,101}

\colorlet{onOffColor}{TUMgreen!90!white}
\colorlet{onOffText}{onOffColor!70!black}
\colorlet{onOffBack}{onOffColor!30!white}

\colorlet{tempColor}{TUMblue}
\colorlet{tempText}{tempColor!80!black}
\colorlet{tempBColor}{TUMred}
\colorlet{tempBText}{tempBColor!80!black}

\colorlet{heating}{red}
\colorlet{heatingBack}{heating!8!white}
\colorlet{heatingText}{heating!70!black}

\colorlet{loadColor}{TUMblue}
\colorlet{residualAColor}{TUMred}
\colorlet{residualBColor}{TUMgreen}

\colorlet{approxColor}{TUMblue}
\colorlet{approxText}{approxColor!80!black}
\colorlet{approxError}{TUMivory!60}

\colorlet{quantileColor}{TUMblue}
\colorlet{medianColor}{TUMred}
\colorlet{quantileBack}{TUMivory!50!white}

\colorlet{Tbase}{TUMblue!80!black}
\colorlet{Temp}{TUMblue!80!white}
\colorlet{3Bin}{TUMred}
\colorlet{1Bin}{TUMgreen}
\colorlet{1BinS}{TUMorange}

\pgfplotsset{
	myboxplot/.style={
		width=8.4cm,
		height=8cm,
		boxplot/draw direction=y,
		xtick={0.5,1.5,2.5,3.5,4.5,5.5},
		x tick label as interval=true,
		xticklabels={1-Bin,1-Bin*,3-Bin,Temp},
		boxplot/every box/.style={fill=quantileBack,draw=quantileColor,mark=quantileColor,thick},
		boxplot/every whisker/.style={draw=quantileColor,thick},
		boxplot/every median/.style={color=medianColor,very thick}
	}
}

\usepackage{hyperref}
\hypersetup{
	breaklinks,
	hidelinks
}

\begin{document}

% Title and Abstract

\title{\texorpdfstring{Improving Accuracy and Efficiency of\\S\kern-0.5pt tart\kern-1pt-\kern-0.5pt up \kern-1pt Cost \kern-1pt Formulations \kern-0.5pt in \kern-1pt MIP \kern-1.5pt Unit \kern-1pt Commitment\\ by Modeling~Power~Plant~Temperatures}{Improving Accuracy and Efficiency of Start-up Cost Formulations in MIP Unit Commitment by Modeling Power Plant Temperatures}}

\author{Matthias~Silbernagl, Matthias~Huber, and Ren\'{e}~Brandenberg
\thanks{}%
\thanks{The authors are with the Technische Uni\-ver\-si\-t\"at M\"un\-chen, 80933 M\"{u}n\-chen, Germany (e-mail: silbernagl@tum.de, matthias.huber@tum.de, brandenb@ma.tum.de)}%
\thanks{}%
}

\markboth{Silbernagl, Huber, and Brandenberg: Improving Accuracy and Efficiency of Start-up~Cost~Models}{}

\begin{textblock}{0.8}(0.1,0.965)
    \noindent
	\footnotesize
	\centering
Accepted for publication in IEEE Transactions on Power Systems, \copyright2015 IEEE
\end{textblock}

\maketitle

\begin{abstract}
This paper presents an improved mixed-integer model for the Thermal Unit Commitment Problem. By introducing new variables for the temperature of each thermal unit, the off-time-dependent start-up costs are modeled accurately and with a lower integrality gap than state-of-the-art formulations. This new approach significantly improves computational efficiency compared to existing formulations, even if they only model a rough approximation of the start-up costs. Our findings were validated on real-world test cases using CPLEX.
\end{abstract}


\IEEEpeerreviewmaketitle

% Content

\section*{Nomenclature}
\label{section:nomenclature}
\addcontentsline{toc}{section}{Nomenclature}

\newcommand{\T}{\mathcal{T}}
\newcommand{\Tsize}{|\T|}
\newcommand{\I}{\mathcal{I}}
\newcommand{\Nodes}{\mathcal{N}}
\newcommand{\Lines}{\mathcal{M}}
\newvar{NodeI}{1}{\mathcal{I}_{#1}}

\renewcommand{\tau}{l}

\vskip-0.3ex

Indices and Sets\\[-2.8ex]
\begin{IEEEdescription}[\IEEEusemathlabelsep\IEEEsetlabelwidth{$CTcccc$}]
	\item[$t \in \squeeze{\T}$] Time periods, $\T = \discrange{1}{T}$
	\item[$i \in \I$] Generating units
	\item[$\tau \in \mathbb{N}$] Look-back time
	\item[$n \in \Nodes$] Network nodes
	\item[$i \in \NodeI(n)$] Generating units at node~$n$
	\item[$m \in \Lines$] Network transmission lines
\end{IEEEdescription}

\newcommand{\MW}{\text{MW}}
\newcommand{\period}{\text{period}}
\newcommand{\cost}{\text{cost}}
\newcommand{\varIndices}[2]{_{#1}^{#2}}

\newvar{demand}{1}{L^{#1}}
\newvar{maxProd}{1}{P^\text{max}_{#1}}
\newvar{minProd}{1}{P^\text{min}_{#1}}
\newvar{maxEfficiency}{1}{\eta^\text{max}_{#1}}
\newvar{minEfficiency}{1}{\eta^\text{min}_{#1}}
\newvar{deltaEfficiency}{1}{\Delta\eta_{#1}}
\newvar{rampup}{1}{RU\mkern-3mu_{#1}}
\newvar{rampdown}{1}{R\mkern-0.5mu D_{#1}}
\newvar{startupRamp}{1}{SU\mkern-3mu_{#1}}
\newvar{shutdownRamp}{1}{S\mkern-2mu D_{#1}}
\newvar{minup}{1}{U\mkern-2mu T_{#1}}
\newvar{mindown}{1}{D\mkern-2mu T_{#1}}

\newvar{fixedProdCost}{1}{A_{#1}}
\newvar{varProdCost}{1}{B_{#1}}

\newvar{StartupCost}{2}{K\varIndices{#1}{#2}}
\newvar{StartupCostTilde}{2}{K\mathllap{\widetilde{\raisebox{-0.2ex}{\phantom{K}}}}\vphantom{K}\varIndices{#1}{#2}}
\newcommand{\StartupTol}{K_{\text{tol}}}

\newvar{heatloss}{1}{\lambda_{#1}}
\newvar{heatingCost}{1}{V_{#1}}
\newvar{fixedStartupCost}{1}{F_{#1}}

\newvar{preOffline}{1}{P\mkern-2.5mu D_{#1}}

\newvar{demandFactor}{1}{\gamma_{#1}}
\newvar{shiftFactor}{2}{\alpha_{#1}^{#2}}
\newvar{lineCapacity}{1}{F_{#1}}

Parameters\\[-2.8ex]
\begin{IEEEdescription}[\IEEEusemathlabelsep\IEEEsetlabelwidth{$CTcccc$}]
	% Production-related
	\item[$\demand(t)$] Electricity demand [$\MW$]
	\item[$\maxProd(i)$] Maximum power output [$\MW$]
	\item[$\minProd(i)$] Minimum power output [$\MW$]
	\item[$\rampup(i)$] Maximum ramp up speed [$\sfrac{\MW}{\period}$]
	\item[$\rampdown(i)$] Maximum ramp down speed [$\sfrac{\MW}{\period}$]
	\item[$\startupRamp(i)$] Maximum ramp up at start-up [$\MW$]
	\item[$\shutdownRamp(i)$] Maximum ramp down at shutdown [$\MW$]
	\item[$\minup(i)$] Minimum uptime [$\period$]
	\item[$\mindown(i)$] Minimum downtime [$\period$]
	\item[$\varProdCost(i)$] Variable production cost [$\sfrac{\cost}{\MW\period}$]
	\item[$\fixedProdCost(i)$] Fixed cost while online [$\sfrac{\cost}{\period}$]
%
	% Start-up cost related
	\item[$\StartupCost(i,\tau)$] Start-up cost after $\tau$ offline periods [$\cost$]
	\item[$\StartupTol$] Start-up cost approximation tolerance (relative)
	\item[$\heatloss(i)$] Heat-loss coefficient, $\heatloss(i) \in (0,1)$ [$\sfrac{1}{\period}$]
	\item[$\heatingCost(i)$] Variable start-up cost [$\cost$] 
	\item[$\fixedStartupCost(i)$] Fixed start-up cost [$\cost$]
	\item[$\preOffline(i)$] Offline periods prior to first period [$\period$]

	% Network related
	\item[$\demandFactor(n)$] Relative demand in nodes
	\item[$\shiftFactor(m,n)$] Relative flow on line~$m$ due to export in node~$n$
	\item[$\lineCapacity(m)$] Transmission line capacity [$\MW$]
\end{IEEEdescription}

\newvar{onOff}{2}{v\varIndices{#1}{#2}}
\renewvar{prod}{2}{p\varIndices{#1}{#2}}
\newcommand{\tempName}{\text{temp}}
\newvar{temp}{2}{\tempName\varIndices{#1}{#2}}
\newvar{tempFunc}{1}{\tempName_{#1}}
\newvar{tempDisc}{2}{\widehat{\tempName}\vphantom{p}\varIndices{#1}{#2}}
\newvar{heating}{2}{h\varIndices{#1}{#2}}
\newvar{shutdownIndicator}{2}{z\varIndices{#1}{#2}}
\newvar{startupIndicator}{2}{y\varIndices{#1}{#2}}
\newvar{startupType}{3}{\delta\varIndices{#1}{#2}(#3)}

\newvar{prodCost}{2}{cp\varIndices{#1}{#2}}
\newvar{startupCost}{2}{cu\varIndices{#1}{#2}}

\newvar{onOffT}{2}{\tilde{v}\varIndices{#1}{#2}}
\newvar{tempT}{2}{\tilde{\text{temp}}\varIndices{#1}{#2}}
\newvar{heatingT}{2}{\tilde{h}\varIndices{#1}{#2}}
\newvar{startupCostT}{2}{\tilde{cu}\varIndices{#1}{#2}}

Variables\\[-2.8ex]
\begin{IEEEdescription}[\IEEEusemathlabelsep\IEEEsetlabelwidth{$CTcccc$}]
	\item[$\onOff(i,t)$] State of power plant, $\onOff(i,t) \in \set{0,1}$
	\item[$\prod(i,t)$] Power output [$\MW$]
	\item[$\temp(i,t)$] Temperature, normalized to $[0,1]$
	\item[$\heating(i,t)$] Heating, normalized to $[0,1]$
	\item[$\shutdownIndicator(i,t)$] Shutdown status, $\shutdownIndicator(i,t) \in \set{0,1}$
	\item[$\startupIndicator(i,t)$] Start-up status, $\startupIndicator(i,t) \in \set{0,1}$
	% Not listed, since not part of our model:
%	\item[$\shutdownIndicator(i,t)$] Shutdown status, $\shutdownIndicator(i,t) \in \set{0,1}$
%	\item[$\startupType(i,t,l)$] Start-up type, $\startupType(i,t,\tau) \in \set{0,1}$
	\item[$\prodCost(i,t)$] Production costs [$\cost$]
	\item[$\startupCost(i,t)$] Start-up costs [$\cost$]
\end{IEEEdescription}  % Notation
\section{Introduction}

\IEEEPARstart{R}{enewable} power sources are being introduced in many of the world's power systems \cite{IEA2013}. The intermittent nature of the power production from these sources, especially of wind and solar, is leading to a higher number of start-ups of conventional thermal power plants \cite{Keatley2013,Huber2013}. Consequently, the percentage of costs caused by start-ups is increasing and accurate start-up cost models are gaining importance.

Operational planning of power systems includes the scheduling of power generating units, which is known as the Unit Commitment (UC) problem \cite{Baldick1995}. Finding cost-optimal solutions to this problem has been an active field of research since almost the beginning of electrification, and a wide variety of optimization approaches have been applied \cite{Sheble1994,Padhy2004}.

A prevalent employed approach is Mixed Integer Programming (MIP) by Branch\&Cut, which is known for simultaneously producing a series of improving solutions and reducing their worst-case optimality gap, leading to an optimal solution. Its main drawback is the high computational effort which has been mitigated by new UC formulations, faster solvers, and greater computational power; still, further progress is vital. This paper contributes by improving the formulation of the start-up costs.

\subsection{Literature Review}

A widely used Unit Commitment model was presented in~\cite{Carrion2006}, describing a formulation of the start-up costs based on \cite{Nowak2000}. Since then, numerous advancements have been published. We restrict ourselves to mentioning only those with a focus similar to our work. The start-up types introduced in \cite{Muckstadt1968} are enhanced in \cite{Simoglou2010} to model the start-up process, including synchronization times, soak times, and power trajectories. Even when disregarding the start-up process,  these start-up types lead to tighter formulations of the start-up costs (see~\cite{Morales-Espana_Thermal_2013}).

Tighter UC formulations have been of interest in general. In \cite{Lee2004}, minimal up- and down-time constraints are considered, proving that the feasible operational schedules can be described by $\bigO(2^T|\I|)$ inequalities. By using start-up and shutdown status variables, \cite{Rajan2005} characterizes the same feasible set with $\bigO(|\I|T)$ inequalities---an example of how representations of polytopes may be simplified by introducing additional variables.

The quadratic production costs have commonly been modeled by piece-wise linear approximations. \cite{Frangioni2009} presents tight valid inequalities for such costs, enabling an iterative approximation scheme. A similar approximation scheme with different valid inequalities is given in \cite{Viana2013}.
Finally, in \cite{Ostrowski2012} solution times are improved by using valid inequalities for the ramping process.

\subsection{Contribution and Paper Organization}

The focus of this paper is a novel approach to modeling the start-up costs. After a short recapitulation of the prevalent state-of-the-art formulations in Section~\ref{section:ExistingModel}, our contributions are introduced in the following order.

Section~\ref{section:tightening} introduces a simple modification of the start-up cost model as presented in \cite{Nowak2000}, which reduces the integrality gap of the model.

While current UC formulations are capable of accurately modeling any increasing start-up cost function, the start-up costs are generally approximated by a step function to keep computation times reasonable (2-5 steps in \cite{Nowak2000,Carrion2006,Simoglou2010,Ostrowski2012,Morales-Espana_Thermal_2013,Morales-Espana_Startup_2013}). Section~\ref{section:approximation} describes an algorithm which chooses a minimal set of steps approximating a given function to a given tolerance (previously published in \cite{whitepaper}).

Section~\ref{section:ThermalStartupCosts} examines the derivation of the commonly used start-up cost function from a simple physical model. This motivates the formulation presented in the following section and serves as a satisfying interpretation.

Section~\ref{section:TemperatureModel} presents the new approach which explicitly models the cooling behavior of units during the offline time by introducing temperature variables. After shutting down a unit, its temperature decays exponentially. At start-up, the lost thermal energy, which is proportional to the temperature loss, must be compensated for by burning additional fuel. By internalizing this physical process instead of encapsulating it in the start-up cost parameters $\StartupCost(i,l)$, this formulation is able to model exact start-up costs for arbitrarily long offline times. Moreover, it significantly improves computational performance compared to existing formulations by  considerably reducing the integrality gap.

Section~\ref{section:NumericalExamples} lists results of numerical experiments, which clearly show the advantages of the proposed approach.  % Introduction
\section{State of the Art}
\label{section:ExistingModel}

This section describes the two prevalent MIP models in recent publications: the approach of \cite{Nowak2000,Carrion2006} with 1 binary variable per unit and period (\enquote{1-Bin}) as well as the approach with 3 binaries (\enquote{3-Bin}) according to \cite{Muckstadt1968,Simoglou2010}, which proved to model start-up costs more efficiently \cite{Morales-Espana_Thermal_2013} and may be extended to model the start-up process.

\subsection{Base Model}
\label{section:CostsAndConstraints}

Each of the discussed start-up cost models may be embedded in any Unit Commitment formulation which represents the operational state of unit~$i$ in period~$t$ with a binary variable~$\onOff(i,t)$, and which minimizes the start-up costs. In the following, we give two UC formulations that are used for comparing the impact of the start-up cost model on the computational performance in Section~\ref{section:NumericalExamples},
\begin{itemize}
	\item a basic formulation without start-up/shutdown indicators, and
	\item an extended formulation which adds indicators, tighter ramping, minimum up-/downtime, and transmission lines.
\end{itemize}

Both formulations share the goal of fulfilling the electricity demand~$\demand(t)$ at minimal cost, which comprises production costs~$\prodCost(i,t)$ and start-up costs~$\startupCost(i,t)$. Denoting the production of each unit as $\prod(i,t)$, this may be modeled as
\newlength{\helper}
\setlength{\helper}{(\widthof{$\scriptstyle i \in \I, t \in \T$}-\widthof{$\sum$})/2}
\begin{align}
	\label{equation:ObjectiveFunction}
	\min \hspace{\helper} &\sum_{\mathclap{i \in \I, t \in \T}} \hspace{\helper} \prodCost(i,t) + \startupCost(i,t) \ \ \text{s.t.}\\
	\label{equation:Demand}
	&\sum_{\mathclap{i \in \I}} \prod(i,t) = \demand(t) \qquad \forall\,t \in \T.
\end{align}

The start-up costs are discussed in Section~\ref{section:StartupCosts}. We use the production costs in \cite{Morales-Espana_Thermal_2013}, which depend linearly on the binary operational state~$\onOff(i,t)$ and the production~$\prod(i,t)$:
\begin{equation}
	\label{constraint:ProductionCost}
	\prodCost(i,t) = \fixedProdCost(i) \onOff(i,t) + \varProdCost(i) \prod(i,t) \qquad \forall\,i \in \I, t \in \T.
\end{equation}

Generally used constraints of power plants regard the minimal production~$\minProd(i)$, the maximal production~$\maxProd(i)$, maximal up and down ramping speeds $\rampup(i)$ and $\rampdown(i)$, the maximal ramping at start-up $\startupRamp(i)$ and shutdown $\shutdownRamp(i)$,  as well as the minimum uptime~$\minup(i)$ and downtime~$\mindown(i)$. The production limits are formulated as
\begin{equation}
	\label{constraint:ProdLimits}
	\minProd(i) \onOff(i,t) \leq \prod(i,t) \leq \maxProd(i) \onOff(i,t) \qquad\quad \forall\,i \in \I, t \in \discrange{1}{T}.
\end{equation}

In the basic formulation without start-up/shutdown indicators, the ramping is modeled as in \cite{Carrion2006} (\cf Appendix~\ref{section:Appendix}, \eqref{constraint:RampUp}--\eqref{constraint:Shutdown}).

The extended formulation uses the tighter ramping due to \cite{Ostrowski2012} (\cf Appendix~\ref{section:Appendix}, \eqref{equation:OstrowskiStart}--\eqref{equation:OstrowskiEnd}) and the tight description of the minimum up-/downtime due to \cite{Rajan2005} (\cf Appendix~\ref{section:Appendix}, \eqref{equation:minup1},\eqref{equation:minup2}). To model the transmission network, the demand is distributed to the nodes~$n$ by the factors~$\demandFactor(n)$, and the net export at each node~$n$ is distributed to the lines~$m$ by the power transfer distribution factors~$\shiftFactor(m,n)$ (see \eg \cite{Bergh2014}). It then suffices to enforce the line capacities~$\lineCapacity(m)$ by
\begin{IEEEeqnarray}{lr}
	\IEEEeqnarraymulticol{2}{c}{
		\label{constraint:transmission}
		\hspace*{3em}- \lineCapacity(m) \leq \mkern-4mu\sum_{n \in \Nodes} \mkern-6mu \shiftFactor(m,n) \bigg(\sum_{i \in \NodeI(n)(n)} \mkern-10mu \prod(i,t) - \demandFactor(n)\demand(t)\bigg) \leq \lineCapacity(m)\hspace*{3em}
	}\\
	&\forall\, t \in \T, m \in \Lines.\nonumber
\end{IEEEeqnarray}
\vspace*{-6mm}

\subsection{Start-up Costs}
\label{section:StartupCosts}

The start-up costs depend on the amount of time~$\tau$ that a unit has been offline before a start-up. For thermal units, they are typically  modeled according to \eg \cite[p. 154]{Wood2013},\cite{Spliethoff2010} as
\begin{equation} \label{equation:StartupCost}
	\StartupCost(i,\tau) = \heatingCost(i)(1 - \exp{-\heatloss(i)\tau}) + \fixedStartupCost(i) \qquad \forall\,i\in \I,\tau \in \mathbb{N},
\end{equation}
where $\fixedStartupCost(i)$ are the fixed start-up costs and $\heatingCost(i)$ are the maximum variable start-up costs, such that the costs for a complete cold start are $\heatingCost(i)+\fixedStartupCost(i)$. The fixed costs include labor costs as well as time-independent wear and tear costs. As the modeled time range is discretized into periods, only integer offline times $l \in \mathbb{N}$ may occur (\cf Fig.~\ref{figure:Startupcosts}).

\subsubsection{Formulation with one binary variable (``1-Bin'')}
\label{section:StartupCosts1bin}

The cost function is modeled by an increasing step function, \ie a piece-wise constant, increasing function. According to \cite{Carrion2006} and \cite{Nowak2000}, this can be formulated as:
\begin{IEEEeqnarray}{c}
	\label{constraint:startcostorig}
	\startupCost(i,t) \geq \StartupCost(i,\tau) \Big(\onOff(i,t) - \mkern-3mu\sum_{n=1}^{\tau} \onOff(i,t-n) \Big)\\
	\IEEEeqnarraymulticol{1}{r}{
		\IEEEnonumber
		\mkern106mu\forall\,i \in \I,t \in \T,\tau \in \discrange*{1}{t-1} \text{ with } \StartupCost(i,l) > \StartupCost(i,l-1)\mkern-3mu.}
\end{IEEEeqnarray}

\subsubsection{Formulation with three binary variables (``3-Bin'')}
\label{section:StartupCosts3bin}

The authors of \cite{Simoglou2010} and \cite{Morales-Espana_Thermal_2013} show that by using the start-up status~$\startupIndicator(i,t)$ and shutdown status~$\shutdownIndicator(i,t)$ as described in \cite{Garver1962},
\begin{IEEEeqnarray}{lll}
	\label{constraint:logic1}
	\startupIndicator(i,1) - \shutdownIndicator(i,1) = \begin{cases}
		\onOff(i,1) & \text{if $\preOffline(i) > 0$,}\\
		\onOff(i,1) - 1 & \text{else,}
	\end{cases} &&\forall\,i \in \I,\\
	\label{constraint:logic2}
	 \startupIndicator(i,t) - \shutdownIndicator(i,t) = \onOff(i,t) - \onOff(i,t-1) &\quad&\forall\,i \in \I, t \in \discrange*{2}{T},\IEEEeqnarraynumspace
\end{IEEEeqnarray}
the start-up costs may be modeled computationally more efficient than using solely $\onOff(i,t)$ as in 1-Bin. To this end, for each unit~$i$ the off-times~$\discrange*{1}{T-1}$ are grouped by their start-up costs into a minimal number~$S_i$ of intervals $L_i^1 \mkern1mu\dot\cup \ldots \dot\cup\mkern2mu L_i^{S_i} = \discrange*{1}{T-1}$ with 
\begin{equation*}
	\forall\,i \in \I, s \in \discrange*{1}{S_i}, \tau, \tau' \in L_i^s: \qquad \StartupCost(i,\tau) = \StartupCost(i,\tau').
\end{equation*}
If unit~$i$ starts up in period~$t$ after $l$~offline periods with $l \in L_i^s$, then the start-up has type~$s$, expressed by $\startupType(i,t,s) = 1$. According to \cite{Simoglou2010} this is modeled by
\begin{IEEEeqnarray}{rClll}
	\label{constraint:onestartup}
	\startupIndicator(i,t) &=& \sum_{s \in \discrange{1}{S_i}} \mkern-11mu \startupType(i,t,s) &&\forall\,i \in \I,t \in \T,\\
	\label{constraint:StartupType}
	\startupType(i,t,s) &\leq& \sum_{\tau \in L_i^s} \shutdownIndicator(i,t-\tau) &\qquad&\forall
		\begin{aligned}[t]
			&\,i \in \I, s \in \discrange*{1}{S_i-1}, \\
			&t \in \T \text{ with } t > \max L^s_i.
		\end{aligned}
	\IEEEeqnarraynumspace
\end{IEEEeqnarray}
While $\startupType(i,t,s)$ may be used to model the start-up process \cite{Simoglou2010}, this comparison considers only the start-up costs by substituting the variables~$\startupCost(i,t)$ in the objective function~\eqref{equation:ObjectiveFunction} by 
\begin{IEEEeqnarray}{lll}
	\label{constraint:startcosttype}
	\startupCost(i,t) := \mkern-4mu\sum_{s \in \discrange{1}{S_i}} \StartupCost(i,\min\!L_i^s) \startupType(i,t,s) \qquad \forall\,i \in \I,t \in \T.
\end{IEEEeqnarray}

\begin{figure}[b]
	\centering
		\begin{tikzpicture}[x=1cm, y=1cm]
		% Gitter zeichnen
		\def\F{0.8}
		\def\V{3.6}
		\def\fc(#1){(\F+\V*(1-exp(-(#1)*2.7/\T)))}
		\def\T{30}
		\def\xscale{7.2/\T}

		\def\a{5}
		\def\b{12}
		\def\c{\T}
		\def\p{{(0.5*\fc( 0) + 0.5*\fc(\a-1))}}
		\def\q{{(0.5*\fc(\a) + 0.5*\fc(\b-1))}}
		\def\r{{(0.5*\fc(\b) + 0.5*\fc(\c-1))}}
		
		\begin{scope}[xscale=\xscale]
			% Fehler zeichnen
			\fill[approxError] (0,\F) \foreach \x in {1,2,...,\T} {
				|- (\x,{\fc(\x-1)})
			} |- (\b,\r) |- (\a,\q) |- (0,\p) -- cycle;

			% Originale Kosten
			\foreach \x in {1,2,...,\T} {
				\draw[very thick,line cap=butt] (\x-1,{\fc(\x-1)}) -- (\x,{\fc(\x-1)});
			}
			
			% Approximation
			\draw[approxColor,very thick,line cap=butt]
				(0,\p) -- (\a,\p)
			    (\a,\q) -- (\b,\q)
			    (\b,\r) -- (\c,\r);
			
			% Stufen
			\begin{scope}[yshift=-0.5cm,approxText,anchor=south]
				\node at ({\a/2},0) {hot};
				\node at ({(\a+\b)/2},0) {warm};
				\node at ({(\b+\c)/2},0) {cold};
			\end{scope}
			\foreach \x in {\a,\b,\c} {
				\draw[thick] (\x,0.13) -- (\x,-0.13);
			}
			\foreach \x in {1,2,...,\T} {
				\ifthenelse{\not\(\x=\a \OR \x=\b \OR \x=\c\)}{
					\draw (\x,0.05) -- (\x,-0.05);
				}{}
			}

			% Funktionen beschriften
			\node[anchor=south east] at (\T,{\fc(\T-1)}) {$\StartupCost(i,l)$};
			\node[anchor=north east,approxText] at (\T,\r) {$\StartupCostTilde(i,l)$};
		\end{scope}

		% Approximationsfehler
		\draw [decorate,xshift=0.8mm,decoration={brace,amplitude=1mm,mirror}]  ({\a*\xscale},\p) -- ({\a*\xscale},{\fc(\a-1)}) node[midway,anchor=west,xshift=0.75mm] {approximation error};
		
		% Achsen zeichnen
		\draw[thick,->] (-0.3,0) -- (7.7,0);
		\begin{pgfinterruptboundingbox}
			\node[anchor=west] at (7.7,0) {$\tau$};
		\end{pgfinterruptboundingbox}
		\draw[thick,->] (0,-0.3) -- (0,{\fc(\T-1)+0.2}) node[anchor=south] {cost}; 
	\end{tikzpicture}
	\caption[Exact discrete start-up costs and three-step approximation.]{Exact discrete start-up costs~$\StartupCost(i,\tau)$ and a three-step approximation~$\StartupCostTilde(i,\tau)$.}
	\label{figure:Startupcosts}
\end{figure}

\section{Improving the 1-Bin and 3-Bin Formulations}
\label{section:ImprovingStepFunction}

In this section, we present a modification of the constraints in (\ref{constraint:startcostorig}) that tightens 1-Bin, and a method to control the approximation error of the time-dependent start-up costs in both 1-Bin and 3-Bin.

\subsection{Tightening the 1-Bin Formulation}
\label{section:tightening} 

The inequalities~(\ref{constraint:startcostorig}) can be tightened by lessening the coefficients of the variables~$\onOff(i,t-n)$,
\begin{IEEEeqnarray}{c}
	\label{constraint:startcosttight}
	\startupCost(i,t) \geq \StartupCost(i,\tau) \onOff(i,t) - \sum_{n=1}^{\tau} (\StartupCost(i,\tau) - \StartupCost(i,n-1)) \onOff(i,t-n)\\
	\IEEEeqnarraymulticol{1}{r}{
		\IEEEnonumber
		\mkern106mu\forall\,i \in \I,t \in \T,\tau \in \discrange*{1}{t-1} \text{ with } \StartupCost(i,l) > \StartupCost(i,l-1)\mkern-3mu.}
\end{IEEEeqnarray}

Each of these inequalities is trivially fulfilled if unit~$i$ is offline in period~$t$, since then its right-hand side is non-positive. If unit~$i$ is online in period~$t$, consider all $n \in \discrange{1}{\tau}$ with $\onOff(i,t-n) = 1$. If no such~$n$ exists, then both the start-up costs~$\startupCost(i,t)$ and the right-hand side of the inequality equal $\StartupCost(i,\tau)$. Otherwise, choose a minimal $n$ with this property. Then, the start-up costs~$\startupCost(i,t)$ equal $\StartupCost(i,n-1)$ and the right-hand side is at most $\StartupCost(i,n-1)$. As these inequalities dominate those in (\ref{constraint:startcostorig}), \ie as they provide a stronger bound on $\startupCost(i,t)$, they still properly model the start-up costs.

The impact of the tightening on the integrality gap is depicted in Fig.~\ref{figure:IntGapRel} in Section~\ref{section:NumericalExamples}.

\subsection{Approximating the Time-dependent Start-up Costs}
\label{section:approximation}

To keep computational efforts reasonable, the time-dependent start-up costs are often approximated either by a constant value (see \eg \cite{Garver1962}) or by up to three steps, \mbox{hot-,} \mbox{cold-,} and possibly warm-start (see \eg \cite{Simoglou2010}). In the light of cool-down times of up to 120~hours for large thermal units \cite{Spliethoff2010}, these approaches result in considerable approximation errors. This is addressed in \cite{Muckstadt1968,Nowak2000} and subsequently in \cite{Carrion2006,Simoglou2010,Ostrowski2012,Morales-Espana_Thermal_2013}, where an arbitrary number of steps is allowed.

When solving MIPs, the goal typically is to reach a certain maximal relative optimality gap. Hence, the approximation $\StartupCostTilde(i,\tau)$ of $\StartupCost(i,\tau)$ needs to guarantee a maximal relative error $\StartupTol$,
\begin{equation}
	\label{equation:accuracy}
	\left| \StartupCostTilde(i,\tau) - \StartupCost(i,\tau) \right| \leq \StartupTol \cdot \StartupCost(i,\tau) \qquad \forall\,i \in \I,\tau \in \discrange*{1}{T-1}.
\end{equation}

In \cite{whitepaper}, we present an algorithm which determines, given a certain approximation tolerance~$\StartupTol$, how to choose $\StartupCostTilde(i,\tau)$ with a minimal number of steps.  % Current model, tightening and thinning
\section{Start-up Costs of Thermal Units}
\label{section:ThermalStartupCosts}

The step-wise start-up cost models considered in the previous section are applicable for all increasing start-up cost functions. However, as is mentioned in the last section, the start-up cost function of a thermal unit is commonly (see \eg \cite[p. 154]{Wood2013},\cite{Spliethoff2010}) defined much more restrictively as
\begin{equation}
	\label{equation:StartupCostFunction}
	\StartupCost(i,\tau) = \underbrace{\heatingCost(i)\,(1-\exp{-\heatloss(i) \tau})}_{\text{variable cost}} + \underbrace{\fixedStartupCost(i)}_{\text{fixed cost}} \qquad \forall\,i \in \I, \tau \in \mathbb{N},
\end{equation}
where $\tau$ denotes the offline time.

The constant costs are derived from the start-up status~$\startupIndicator(i,t)$ modeled as in constraints~\eqref{constraint:logic1}, \eqref{constraint:logic2}. The variable costs originate from the reheating process at start-up, where fuel needs to be burned and where the unit experiences thermal stress.

Here, the term~$(1 - \exp{-\heatloss(i) \tau})$ is proportional to the heat loss of the power plant incurred while offline, and models the exponential decay of the temperature,
\begin{equation}
	\label{equation:ContinuousTemperature}
	\tempFunc(i)(\tau) = \exp{-\heatloss(i) \tau} \qquad \forall\,i \in \I, \tau \in \mathbb{R}_{\geq 0},
\end{equation}
assuming the operational temperature is normalized to~$1$ and the environmental temperature is normalized to~$0$.

\begin{figure}[t]
	\centering
	\begin{tikzpicture}[x=0.7cm,y=0.7cm]
		\def\T{10}
		\def\lambda{0.32}
		\def\h{3}
		\def\fT(#1){(exp(-\lambda*(#1))}
		\def\fdT(#1){\lambda*(exp(-\lambda*(#1)))}
		\def\ftemp(#1){((#1) <= 1 ? 1 : ((#1) < 4 ? \fT(#1-1) : ((#1) <= 5 ? 1 : \fT(#1-5))))}
		
		\begin{scope}[yshift=1cm]
			\begin{scope}[yscale=\h]
				% continuous temperature (samples=10 prevents error)
				\begin{scope}[tempColor,very thick, line cap=round, line join=round]
					\draw (1,1) plot[smooth, domain=1:4, samples=10]  (\x, {\fT(\x-1)});
					\draw[loosely dashed, dash phase=1.1pt] (4,{\fT(3)}) -- (4,1);
					\draw (5,1) plot[smooth, domain=5:10, samples=10] (\x, {\fT(\x-5)});
				\end{scope}
	
				% discrete temperature
				\foreach \k in {1,...,10} {
					\coordinate (temp\k) at ({\k-1},{\ftemp(\k-1)});
					%\fill[tempColor] (temp\k) ellipse (0.07 and {0.07/\h});
					\draw[tempBColor, very thick, line cap=round] (temp\k) -- +(1,0);				
				}
				\draw (temp7) ++(0.4,0.05) node[tempBText,anchor=south west] (tempDisc) {$\tempDisc(i,t)$};
				\draw (tempDisc |- 0,0.2) node[tempText] {$\tempFunc(i)(t)$};
	
				% heating
				\draw[tempBColor, very thick, line cap=round, loosely dashed, dash phase=0pt] (4, {\fT(3)}) -- (5, {\fT(3)});
				\draw[tempBColor, very thick, -latex] (4.53, {\fT(3)}) -- (4.53, 1);
				\node[tempBText, anchor=west] at (4.58, {(\fT(3)+1)/2}) {$\heating(i,4)$};

				% tick for 1
				\draw (0.08, 1) -- (-0.08, 1) node[anchor=east] {$1$};
				\coordinate (max) at (0,1);
			\end{scope}
	
			% first x axis
			\draw[<->] (-0.5, 0) -- (\T + 0.5, 0);
			\foreach \i in {1,...,\T} {
				\draw ({max(0,\i)+min(0,\i-1)}, 0) +(0, -0.08) -- +(0, 0.08);
				    ++( 0, 2) +(0, -0.08) -- +(0, 0.08);
			}
		\end{scope}
		
		\begin{scope}[yscale=1]
			% onOff
			\draw[onOffColor, very thick,fill=onOffBack] (0,0) -- (0,1) -| (1,0) -| (4,1) -| (5,0) -- (\T,0);
			\node[onOffText, anchor=west] at (1, 0.5) {$\onOff(i,)$};

			% tick for 1
			\draw (0.08, 1) -- (-0.08, 1) node[anchor=east] {$1$};
		\end{scope}
		
		% second x axis
		\draw[<->] (-0.5, 0) -- (\T + 0.5, 0);
		%\node[anchor=west] at (\T + 0.5, 0) {$t$};
		%\node[anchor=north east] at (0,0) {$0$};
		\node[anchor=west] at (\T + 0.5, 0) {$t$};
		\foreach \i in {1,...,\T} {
			%\node[anchor=north] at (\i, -0.08) {$\i$};
			\node[anchor=north] at ({\i-0.5}, -0.08) {$\i$};
			\draw ({max(0,\i)+min(0,\i-1)}, 0) +(0, -0.08) -- +(0, 0.08);
			    ++( 0, 2) +(0, -0.08) -- +(0, 0.08);
		}
		
		% y axis
		\draw[<-] (max) ++(0,0.5) -- (0, -0.5);
	\end{tikzpicture}

	\caption[Discretization of a unit's temperature function.]{Discretization of a unit's temperature function. Following the operational schedule, the unit exhibits the temperature function~$\tempFunc(i)$ which is discretized to $\tempDisc(i,)$, with resulting heating~$\heating(i,)$ according to (\ref{constraint:TemperatureDevelopmentFirst}).}
	\label{figure:Temperature}
\end{figure}
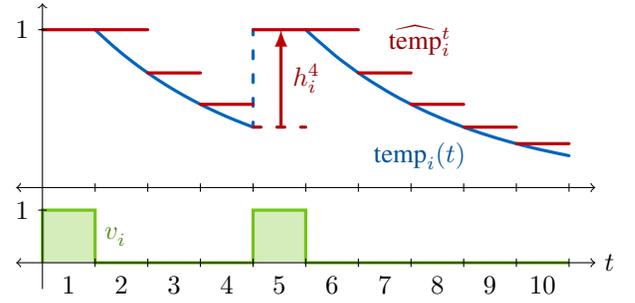

As shown in Fig.~\ref{figure:Temperature}, (\ref{equation:ContinuousTemperature}) is discretized by a step-wise constant function with steps according to
\begin{IEEEeqnarray}{rCl}
	\label{equation:DiscretizedTemperature}
	\tempDisc(i,t) &:=& \begin{cases}
					1 & \text{if $\onOff(i,t) = 1$,}\\
					\tempFunc(i)(\tau\varIndices{i}{t}) = \exp{-\heatloss(i) \tau_i^t} & \text{else,}
	\end{cases}\\
	\IEEEeqnarraymulticol{3}{r}{\forall\,i \in \I,t \in \T,\IEEEnonumber}
\end{IEEEeqnarray}
where $\tau\varIndices{i}{t}$ denotes the number of periods that unit~$i$ is offline prior to period~$t$. 

The above nonlinear definition of $\tempDisc(i,t)$ may be restated recursively as
\begin{IEEEeqnarray}{rCl}
	\label{equation:TemperatureRecursionStart}
	\tempDisc(i,1) &=& \begin{cases}
		1 & \text{if $\onOff(i,1) = 1$,}\\
		\exp{-\heatloss(i)\preOffline(i)}\hspace*{1.8em} & \text{else,}
	\end{cases} \qquad \forall\,i \in \I,\\
	\label{equation:TemperatureRecursion}
	\tempDisc(i,t) &=& \begin{cases}
		1 & \text{if $\onOff(i,t-1) = 1$ or $\onOff(i,t) = 1$,}\\
		\exp{-\heatloss(i)} \tempDisc(i,t-1) & \text{else,}
	\end{cases}\\
	\IEEEeqnarraymulticol{3}{r}{\forall\,i \in \I,t \in \discrange{2}{T}.\IEEEnonumber}
\end{IEEEeqnarray}  % Thermal start-up costs
\section{The Temperature Model}
\label{section:TemperatureModel}

In this section, we model the temperature loss derived in the last section by explicitly capturing the temperature of a unit as the new state variable~$\temp(i,t)$ and the amount of heating as the new variable~$\heating(i,t)$. Combined with the start-up status~$\startupIndicator(i,t)$ (cf. \ref{section:ThermalStartupCosts}), they are used to model the start-up cost function as defined in equation~(\ref{equation:StartupCost}).

The new variables are continuous and non-negative,
\begin{align}
	\label{constraint:TemperatureNonnegative}
	\temp(i,t) &\in \mathbb{R}_{\geq 0} \qquad \forall\,i \in \I,t \in \T,\\
	\label{constraint:HeatingNonnegative}
	\heating(i,t) &\in \mathbb{R}_{\geq 0} \qquad \forall\,i \in \I, t \in \discrange*{0}{T-1}.
\end{align}
The operational temperature is expressed as
\begin{equation}
	\label{constraint:TemperatureNormalization}
	\onOff(i,t) \leq \temp(i,t) \qquad \forall\,i \in \I,t \in \T,
\end{equation}
enforcing a temperature of at least $1$ during operation. The recursion in equations~(\ref{equation:TemperatureRecursionStart}) and (\ref{equation:TemperatureRecursion}) is modeled as
\begin{IEEEeqnarray}{rCl} 
	\label{constraint:TemperatureDevelopmentFirst}
	\temp(i,1) &=& \exp{-\heatloss(i) \preOffline(i)} + \heating(i,0) \qquad \forall\,i \in \I,\\
	\label{constraint:TemperatureDevelopment}
	\temp(i,t) &=& \exp{-\heatloss(i)}\temp(i,t-1) + (1 - \exp{-\heatloss(i)})\onOff(i,t-1) + \heating(i,t-1) \IEEEeqnarraynumspace\\
	            \IEEEeqnarraymulticol{3}{r}{\forall\, i \in \I, t \in \discrange{2}{T},  \IEEEnonumber}%%
\end{IEEEeqnarray}
\vskip-1ex
\noindent which causes the temperature
\begin{itemize}
  \item to decay exponentially while the unit is offline ($\onOff(i,t) = 0$),
  \item to stay constant at~$1$ while the unit is online ($\onOff(i,t) = 1$), and
  \item to rise by $\heating(i,t)$ if the unit is heating.
\end{itemize}

\noindent Finally, the start-up costs are defined as
\begin{equation} \label{constraint:StartupCost}
	\startupCost(i,t) := \heatingCost(i) \heating(i,t-1) + \fixedStartupCost(i) \startupIndicator(i,t) \qquad \forall\, i \in \I,t \in \T.
\end{equation}

We proceed by explaining the correctness of this model. As shown in \cite{Garver1962}, the constant part of the start-up costs is modeled correctly by using the start-up status~$\startupIndicator(i,t)$.

The temperature losses increase proportionally with \hbox{$\temp(i,t) - \onOff(i,t)$}. Thus, in a cost-minimal solution, heating is applied such that the temperature is minimal while fulfilling $\temp(i,t) \geq \onOff(i,t)$. This entails two consequences:
\begin{enumerate}
  \item Heating is applied only in the period prior to each start-up. Earlier heating could be postponed until this period, thus saving heating costs.
  \item The amount of heating is exactly such that the temperature reaches~$1$. Excessive heating could either be postponed until the period prior to the next start-up, or be avoided if there is no such start-up.
\end{enumerate}

Therefore, in a cost-minimal solution, the temperature variables~$\temp(i,t)$ match the discretized temperatures~$\tempDisc(i,t)$ as given in equation~(\ref{equation:DiscretizedTemperature}), and the start-up costs~$\startupCost(i,t)$ equal~$0$ if unit~$i$ does not start in period~$t$.

Given a cost-minimal solution, assume that unit~$i$ starts up in period~$t$ after $\tau$~offline periods. By period~$t-1$ the unit has cooled down for $\tau-1$~periods, and in period~$t$, the temperature after start-up has to be~$1$ again,
\begin{equation*}
	\temp(i,t-1) \stackrel{(\ref{equation:DiscretizedTemperature})}{=} \exp{-\heatloss(i)(\tau-1)} \quad \text{and}\quad \temp(i,t) \stackrel{(\ref{constraint:TemperatureNormalization})}{=} 1,
\end{equation*}
Thus, the needed heating, considering the further cooling during period~$t-1$, matches the expected temperature loss,
\begin{align*}
	\heating(i,t-1) \: & \stackrel{\mathclap{(\ref{constraint:TemperatureDevelopment})}}{=}\: \temp(i,t) - \exp{-\heatloss(i)}\temp(i,\tau-1) + (1-\exp{-\heatloss(i)})\smash{\underbrace{\onOff(i,t-1)}_{=0}}\\
		&= 1 - \exp{-\heatloss(i)\tau}.
\end{align*}
This means the variable part of the start-up costs is modeled correctly too, leading to $\startupCost(i,t) = \StartupCost(i,\tau)$. Note that in solutions that are not cost-optimal, heating may occur in periods not directly prior to a start-up. Thereby, the start-up costs may be incorrect, as is the case for all other formulations too.

While this model uses new additional variables, it reduces the number of constraints in comparison to 1-Bin, 1-Bin* and 3-Bin, even with a start-up cost approximation tolerance of~$\StartupTol=5\%$ (see Section~\ref{section:Sizes}). Fig.~\ref{figure:IntGapRel} suggests that the integrality gap of this model is on average smaller, while the solution times of the linear relaxation remain comparable to 1-Bin and 3-Bin (see Fig.~\ref{figure:IntGapTime}). Both factors are crucial for the improved number of solved instances shown in Fig.~\ref{figure:scalingcomparesolved}.  % Temperature model
\newcommand{\scalingHeight}{5cm}

\newboolean{IntGapVert}
\IntGapVertfalse
\newboolean{IntGapTimeVert}
\IntGapTimeVertfalse

\section{Numerical Examples}
\label{section:NumericalExamples}

This section presents results from numerical examples which show the benefits of our modeling approach. After introducing the modeling setup and reporting problem sizes and solution times of the linear relaxation, its reduced integrality gap is highlighted. This advantage leads to an overall faster optimization procedure and enables larger models to be solved.

The experiments are performed using the CPLEX solver.

\vspace*{-1mm}
\subsection{Scenarios and Data Description}
\label{section:Scenarios}

We investigate two scenarios, one based on the German power system with 223~units, and one based on the IEEE~118 bus system with 54~units and 118~nodes in a transmission network.

\subsubsection{German power system}

The raising requirements for fossil-fuel power plants, which stem from a more volatile residual load, include more start-ups and hence result in a higher ratio of start-up to operational costs \cite{Keatley2013}. We expect the higher percentage of start-up costs to lead to higher solution times and to increase the advantages of our approach. To consider the impact of a more volatile residual load in the numerical experiments, a forecast scenario for the year 2025 is employed.

We use power plant data based on the German power system of 2014 as published by the German Federal Network Agency in \cite{bnetza2014}, comprising 228 individually controlled power plants. The data is augmented by assumptions regarding minimal production, efficiency, and start-up costs, which are partly based on \cite{Kumar2012,Eurelectric2003,Egerer2014}. As we model the year 2025, all nuclear power plants are phased out in favor of four additional combined cycle gas turbines, reducing the number of plants to 223. \textsuperscript{1} 

The main benefits of this dataset are
\begin{itemize}
  \item an adequate number of power plants, representing the diversity of a real power system, and
  \item detailed thermal start-up cost functions given by coefficients~$\fixedStartupCost(i)$, $\heatingCost(i)$, and $\heatloss(i)$ (see~\eqref{equation:StartupCostFunction}).
\end{itemize}

In addition to power plant data, the model requires data of the residual load, \ie of the difference between load and electricity production from must-run renewable power sources. The load data is taken from ENTSO-E \cite{ENTSOE_Consumption_2007} and scaled to a yearly electricity consumption of 520~TWh. Wind and solar electricity generation profiles are computed based on the NASA~MERRA~database~\cite{MERRA} for the same base year. Afterwards, these profiles are scaled according to the respective installed capacity (50~GW wind, 50~GW solar). Biomass and hydro power plants are assumed to produce at full capacity (5.5~GW biomass, 4.5~GW hydro).

Each experiment is performed using 14 time ranges of length~$\squeeze{T=72}$ (Sections~\ref{section:Sizes},\ref{section:IntegralityGap}), length~$\squeeze{T=240}$ (Section~\ref{section:LinearRelaxation}) or varying length (Section~\ref{section:Scaling}), starting in the $S$-th hour of the year with $S \in  \{624k + 433: k \in \discrange{0}{13}\}$. This set is chosen such that each time range starts at midnight, the time ranges are uniformly distributed over the year 2025, and two time ranges start on any day of the week, respectively.

\subsubsection{IEEE 118 bus system}

This scenario is based on the IEEE 118 bus test case published in \cite{morales2014phd}, and again augmented to include the relevant power plant data. \footnote{The complete datasets including power plants, residual loads, and the transmission grid is available as ancillary files at \url{http://arxiv.org/abs/1408.2644}.} Apart from being well-studied, its major benefit is its realistic transmission network.

The test system provides load values for 24 hours and 20 wind scenarios which are concatenated into a residual load of 480~periods. Since the low average wind production of 5.4\% of the load leads to a lower ratio of start-up to operational costs than in the scenario of the German power system, we expect the advantage of our approach to be less pronounced.

Analogous to the German power system, 14 uniformly distributed starting points are given by $\squeeze[0.5]{S \in \{24k + 1: k \in \discrange{0}{13}\}}$.

\subsection{Compared Model Formulations}

We evaluate our approach by comparing it to the state-of-the-art start-up cost formulations introduced in Section~\ref{section:StartupCosts}. We consider
\begin{enumerate}
	\item 1-Bin: Start-up costs modeled by inequality~\eqref{constraint:startcostorig}.
	\item 1-Bin*: Same as 1-Bin, with the tightened start-up cost inequalities~\eqref{constraint:startcosttight} instead of the original inequalities~\eqref{constraint:startcostorig}.
	\item 3-Bin: Same as 1-Bin, except that start-up cost inequalities~(\ref{constraint:startcostorig}) are replaced by the inequalities~\eqref{constraint:logic1}-\eqref{constraint:StartupType}, and the start-up costs are defined as in \eqref{constraint:startcosttype}.
	\item Temp: New approach with explicit modeling of the power plant temperature, as described in Section~\ref{section:TemperatureModel}, including inequalities~\eqref{constraint:TemperatureNonnegative}-\eqref{constraint:TemperatureDevelopment} and start-up costs defined in \eqref{constraint:StartupCost}.
\end{enumerate}

\noindent These formulations are embedded into the two models described in Section~\ref{section:CostsAndConstraints},
\begin{itemize}
  \item the basic formulation composed of \eqref{equation:ObjectiveFunction}--\eqref{constraint:ProdLimits}, \eqref{constraint:RampUp}--\eqref{constraint:Shutdown}, and
  \item the extended formulation composed of \eqref{equation:ObjectiveFunction}--\eqref{constraint:transmission}, \eqref{equation:OstrowskiStart}--\eqref{equation:minup2}.
\end{itemize}
The basic UC problem uses the German power system, while the extended UC problem requires the IEEE 118 bus system.

In 1-Bin, 1-Bin*, and 3-Bin, the start-up costs are approximated with tolerance~$\StartupTol \in \set{0\%,5\%,20\%}$ (see Section~\ref{section:approximation}). Using $\StartupTol = 0\%$, the modeled start-up costs are equal to Temp, resulting in equivalent problems and solutions, which is required when comparing integrality gaps. With $\StartupTol = 20\%$, the start-up cost functions are approximated very roughly with 2.3~steps on average. Finally, as in the presented scenarios start-up costs amount to up to $10\%$ of the total costs, $\StartupTol = 5\%$ is a sensible choice with a maximal error of $0.5\%$ of the objective value.

\ifIntGapTimeVert\else
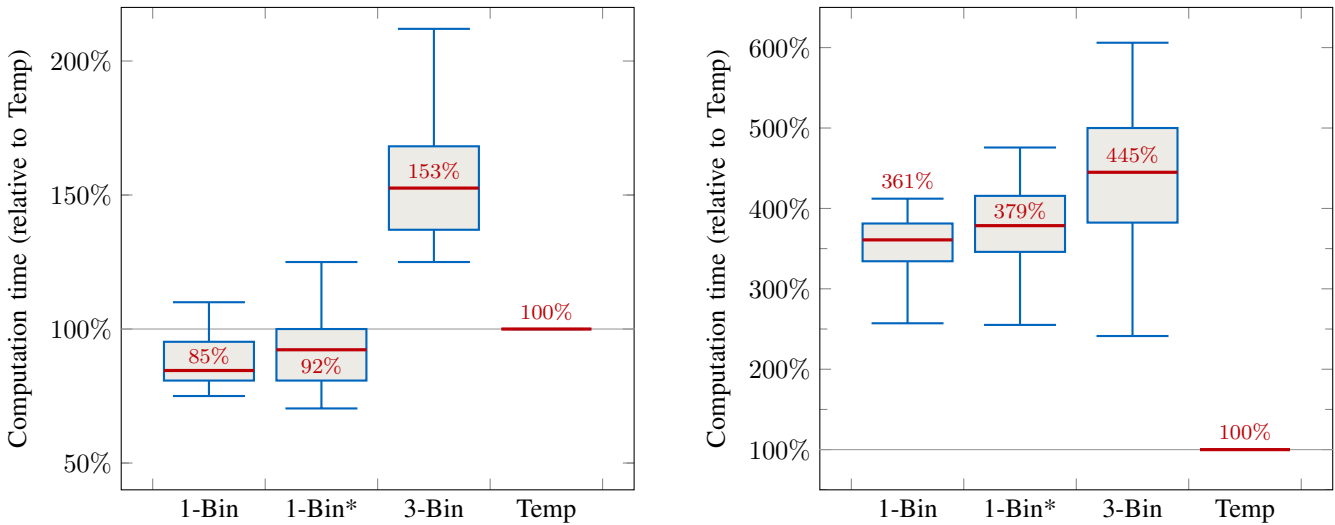
\begin{figure*}[t]
	\hspace*{\fill}
	\begin{tikzpicture}
	\begin{axis}[
		myboxplot,
		ylabel={Computation time (relative to Temp)},
		ymin=0.4,
		ymax=2.2,
		ytick={0.5,1,1.5,2},
		yticklabels={$50\%$, $100\%$, $150\%$, $200\%$},
		minor ytick={0.5,1.5,2.5,3.5}
	]

	\draw[black!40] (axis cs:0,1) -- (axis cs:6,1);

\addplot[mark=+,boxplot prepared={
    lower whisker=0.75,
    lower quartile=0.8076923076923077,
    median=0.8452380952380952,
    upper quartile=0.9523809523809523,
    upper whisker=1.1
}] table[row sep=\\,y index=0] {data\\};
\node[anchor=center,medianColor] at ($(axis cs:1,0.8452380952380952)!0.5!(axis cs:1,0.9523809523809523)$) {\footnotesize$85\%$};
% printBoxplot([0.96 0.8260869565217391 0.8571428571428571 0.8620689655172413               0.76               0.75 0.8076923076923077 0.8636363636363636 0.8148148148148148 0.9523809523809523                  1                1.1 0.8333333333333334 0.7857142857142857])
\addplot[mark=+,boxplot prepared={
    lower whisker=0.7037037037037037,
    lower quartile=0.8076923076923077,
    median=0.9226190476190477,
    upper quartile=1,
    upper whisker=1.25
}] table[row sep=\\,y index=0] {data\\};
\node[anchor=center,medianColor] at ($(axis cs:2,0.9226190476190477)!0.5!(axis cs:2,0.8076923076923077)$) {\footnotesize$92\%$};
% printBoxplot([1                  1 0.9523809523809523 0.8620689655172413                0.8 0.7916666666666666 0.8076923076923077 0.8636363636363636 0.7037037037037037  1.142857142857143                  1               1.25 0.9166666666666666 0.9285714285714286])
\addplot[mark=+,boxplot prepared={
    lower whisker=1.25,
    lower quartile=1.37037037037037,
    median=1.525862068965517,
    upper quartile=1.681818181818182,
    upper whisker=2.12
}] table[row sep=\\,y index=0] {data\\};
\node[anchor=south,medianColor] at (axis cs:3,1.525862068965517) {\footnotesize$153\%$};
% printBoxplot([2.12 1.652173913043478 1.428571428571429 1.551724137931034              1.64 1.291666666666667               1.5 1.409090909090909  1.37037037037037 1.761904761904762 1.681818181818182               1.7              1.25 1.321428571428571])
\addplot[mark=+,boxplot prepared={
    lower whisker=1,
    lower quartile=1,
    median=1,
    upper quartile=1,
    upper whisker=1
}] table[row sep=\\,y index=0] {data\\};
\node[anchor=south,medianColor] at (axis cs:4,1) {\footnotesize$100\%$};
% printBoxplot([1 1 1 1 1 1 1 1 1 1 1 1 1 1])

	\end{axis}
\end{tikzpicture}
	\hspace*{3ex}
	\begin{tikzpicture}
	\begin{axis}[
		myboxplot,
		ylabel={Computation time (relative to Temp)},
		ymin=0.5,
		ymax=6.5,
		ytick={1,2,3,4,5,6},
		yticklabels={$100\%$, $200\%$, $300\%$, $400\%$, $500\%$, $600\%$},
		minor ytick={0.5,1.5,2.5,3.5}
	]

	\draw[black!40] (axis cs:0,1) -- (axis cs:6,1);

\addplot[mark=+,boxplot prepared={
    lower whisker=2.571428571428572,
    lower quartile=3.342857142857143,
    median=3.608823529411765,
    upper quartile=3.8125,
    upper whisker=4.121212121212121
}] table[row sep=\\,y index=0] {data\\};
\node[anchor=south,medianColor] at (axis cs:1,4.1212121212121215) {\footnotesize$361\%$};
% printBoxplot([3.6           3.90625           4.09375 2.571428571428572 3.342857142857143 3.676470588235294 2.760869565217391            3.8125                 3 3.575757575757576 3.638888888888889 4.121212121212121 3.432432432432432 3.617647058823529])
\addplot[mark=+,boxplot prepared={
    lower whisker=2.551020408163265,
    lower quartile=3.45945945945946,
    median=3.785714285714286,
    upper quartile=4.15625,
    upper whisker=4.757575757575758
}] table[row sep=\\,y index=0] {data\\};
\node[anchor=center,medianColor] at ($(axis cs:2,3.785714285714286)!0.5!(axis cs:2,4.15625)$) {\footnotesize$379\%$};
% printBoxplot([3.885714285714286           4.15625           4.40625 2.551020408163265 3.685714285714285 4.058823529411764 2.934782608695652            4.4375 3.512820512820513 4.757575757575758 3.944444444444445 3.424242424242424  3.45945945945946 3.647058823529412])
\addplot[mark=+,boxplot prepared={
    lower whisker=2.41304347826087,
    lower quartile=3.823529411764706,
    median=4.449519230769231,
    upper quartile=5,
    upper whisker=6.060606060606061
}] table[row sep=\\,y index=0] {data\\};
\node[anchor=south,medianColor] at (axis cs:3,4.449519230769231) {\footnotesize$445\%$};
% printBoxplot([3.942857142857143            4.4375             4.375 2.755102040816327 4.885714285714286 3.823529411764706  2.41304347826087             5.875 4.461538461538462 5.090909090909091                 5 6.060606060606061 4.675675675675675 3.764705882352941])
\addplot[mark=+,boxplot prepared={
    lower whisker=1,
    lower quartile=1,
    median=1,
    upper quartile=1,
    upper whisker=1
}] table[row sep=\\,y index=0] {data\\};
\node[anchor=south,medianColor] at (axis cs:4,1) {\footnotesize$100\%$};
% printBoxplot([1 1 1 1 1 1 1 1 1 1 1 1 1 1])

	\end{axis}
\end{tikzpicture}
	\hspace*{\fill}
	\caption{Solution times of the linear relaxation relative to Temp for 14 test cases with $T=120$ periods and $\StartupTol = 5\%$ in German power system (left chart) and the IEEE 118 bus system (right chart). Temp outperforms 3-Bin consistently.}
	\label{figure:IntGapTime}
\end{figure*}
\fi

\subsection{Problem Sizes}
\label{section:Sizes}

Table~\ref{table:Sizes} lists the number of variables and inequalities for the four start-up cost models and for different start-up cost approximation tolerances~$\StartupTol$ in the basic formulation. Their number of additional variables and inequalities remains constant in the extended formulation, except for the start-up and shutdown indicators which are already part of 3-Bin and Temp but have to be added in 1-Bin and 1-Bin*.

While Temp uses twice as many variables as 1-Bin and 1-Bin*, the model requires significantly less inequalities than the state-of-the-art formulations at $\StartupTol = 5\%$. Naturally, higher tolerances~$\StartupTol$ result in fewer inequalities: a number of inequalities approximately equal to Temp is reached by 1-Bin and 1-Bin* at~$\StartupTol \approx 11.3\%$, and by 3-Bin at $\StartupTol \approx 19.2\%$.

\begin{table}[t]
	\vskip-0.3cm
	\caption{Problem sizes for 72~periods and 223~units in basic formulation}
	\label{table:Sizes}
	{
		\centering
		\newcommand{\z}{.00}
\begin{tabular}{lrrrr}
\toprule
Model  & $\StartupTol$ & Avg. steps & Variables & Inequalities\\
\midrule
None   &      &      &  32112 & 79683\\
\addlinespace
1-Bin  &  0\% & 71\z &  48168 & 649671\\
       &  5\% & 6.48 &  48168 & 166334\\
       & 20\% & 2.32 &  48168 & 111423\\
\addlinespace
1-Bin* &  0\% & 71\z &  48168 & 649671\\
       &  5\% & 6.48 &  48168 & 166334\\
       & 20\% & 2.32 &  48168 & 111423\\
\addlinespace
3-Bin  &  0\% & 71\z & 634435 & 665950\\
       &  5\% & 6.48 & 151098 & 182613\\
       & 20\% & 2.32 &  96187 & 127702\\
\addlinespace
Temp &      &      &  96336 & 127851\\
%Temp   &      &      &  96336 & 128818\\
\bottomrule\\
\end{tabular}

	}
	Problem sizes and average number of steps in the approximation of the start-up cost function of all start-up cost formulations for $\StartupTol \in \set{0\%,5\%,20\%}$ for $T = 72$ periods and $223$~units (basic formulation).
\end{table}

\subsection{Computational Effort for Solving the LP}
\label{section:LinearRelaxation}

\ifIntGapTimeVert
\begin{figure}[b]
	\begin{tikzpicture}
	\begin{axis}[
		myboxplot,
		ylabel={Computation time (relative to Temp)},
		ymin=0.4,
		ymax=2.2,
		ytick={0.5,1,1.5,2},
		yticklabels={$50\%$, $100\%$, $150\%$, $200\%$},
		minor ytick={0.5,1.5,2.5,3.5}
	]

	\draw[black!40] (axis cs:0,1) -- (axis cs:6,1);

\addplot[mark=+,boxplot prepared={
    lower whisker=0.75,
    lower quartile=0.8076923076923077,
    median=0.8452380952380952,
    upper quartile=0.9523809523809523,
    upper whisker=1.1
}] table[row sep=\\,y index=0] {data\\};
\node[anchor=center,medianColor] at ($(axis cs:1,0.8452380952380952)!0.5!(axis cs:1,0.9523809523809523)$) {\footnotesize$85\%$};
% printBoxplot([0.96 0.8260869565217391 0.8571428571428571 0.8620689655172413               0.76               0.75 0.8076923076923077 0.8636363636363636 0.8148148148148148 0.9523809523809523                  1                1.1 0.8333333333333334 0.7857142857142857])
\addplot[mark=+,boxplot prepared={
    lower whisker=0.7037037037037037,
    lower quartile=0.8076923076923077,
    median=0.9226190476190477,
    upper quartile=1,
    upper whisker=1.25
}] table[row sep=\\,y index=0] {data\\};
\node[anchor=center,medianColor] at ($(axis cs:2,0.9226190476190477)!0.5!(axis cs:2,0.8076923076923077)$) {\footnotesize$92\%$};
% printBoxplot([1                  1 0.9523809523809523 0.8620689655172413                0.8 0.7916666666666666 0.8076923076923077 0.8636363636363636 0.7037037037037037  1.142857142857143                  1               1.25 0.9166666666666666 0.9285714285714286])
\addplot[mark=+,boxplot prepared={
    lower whisker=1.25,
    lower quartile=1.37037037037037,
    median=1.525862068965517,
    upper quartile=1.681818181818182,
    upper whisker=2.12
}] table[row sep=\\,y index=0] {data\\};
\node[anchor=south,medianColor] at (axis cs:3,1.525862068965517) {\footnotesize$153\%$};
% printBoxplot([2.12 1.652173913043478 1.428571428571429 1.551724137931034              1.64 1.291666666666667               1.5 1.409090909090909  1.37037037037037 1.761904761904762 1.681818181818182               1.7              1.25 1.321428571428571])
\addplot[mark=+,boxplot prepared={
    lower whisker=1,
    lower quartile=1,
    median=1,
    upper quartile=1,
    upper whisker=1
}] table[row sep=\\,y index=0] {data\\};
\node[anchor=south,medianColor] at (axis cs:4,1) {\footnotesize$100\%$};
% printBoxplot([1 1 1 1 1 1 1 1 1 1 1 1 1 1])

	\end{axis}
\end{tikzpicture}
	\vskip1ex
	\begin{tikzpicture}
	\begin{axis}[
		myboxplot,
		ylabel={Computation time (relative to Temp)},
		ymin=0.5,
		ymax=6.5,
		ytick={1,2,3,4,5,6},
		yticklabels={$100\%$, $200\%$, $300\%$, $400\%$, $500\%$, $600\%$},
		minor ytick={0.5,1.5,2.5,3.5}
	]

	\draw[black!40] (axis cs:0,1) -- (axis cs:6,1);

\addplot[mark=+,boxplot prepared={
    lower whisker=2.571428571428572,
    lower quartile=3.342857142857143,
    median=3.608823529411765,
    upper quartile=3.8125,
    upper whisker=4.121212121212121
}] table[row sep=\\,y index=0] {data\\};
\node[anchor=south,medianColor] at (axis cs:1,4.1212121212121215) {\footnotesize$361\%$};
% printBoxplot([3.6           3.90625           4.09375 2.571428571428572 3.342857142857143 3.676470588235294 2.760869565217391            3.8125                 3 3.575757575757576 3.638888888888889 4.121212121212121 3.432432432432432 3.617647058823529])
\addplot[mark=+,boxplot prepared={
    lower whisker=2.551020408163265,
    lower quartile=3.45945945945946,
    median=3.785714285714286,
    upper quartile=4.15625,
    upper whisker=4.757575757575758
}] table[row sep=\\,y index=0] {data\\};
\node[anchor=center,medianColor] at ($(axis cs:2,3.785714285714286)!0.5!(axis cs:2,4.15625)$) {\footnotesize$379\%$};
% printBoxplot([3.885714285714286           4.15625           4.40625 2.551020408163265 3.685714285714285 4.058823529411764 2.934782608695652            4.4375 3.512820512820513 4.757575757575758 3.944444444444445 3.424242424242424  3.45945945945946 3.647058823529412])
\addplot[mark=+,boxplot prepared={
    lower whisker=2.41304347826087,
    lower quartile=3.823529411764706,
    median=4.449519230769231,
    upper quartile=5,
    upper whisker=6.060606060606061
}] table[row sep=\\,y index=0] {data\\};
\node[anchor=south,medianColor] at (axis cs:3,4.449519230769231) {\footnotesize$445\%$};
% printBoxplot([3.942857142857143            4.4375             4.375 2.755102040816327 4.885714285714286 3.823529411764706  2.41304347826087             5.875 4.461538461538462 5.090909090909091                 5 6.060606060606061 4.675675675675675 3.764705882352941])
\addplot[mark=+,boxplot prepared={
    lower whisker=1,
    lower quartile=1,
    median=1,
    upper quartile=1,
    upper whisker=1
}] table[row sep=\\,y index=0] {data\\};
\node[anchor=south,medianColor] at (axis cs:4,1) {\footnotesize$100\%$};
% printBoxplot([1 1 1 1 1 1 1 1 1 1 1 1 1 1])

	\end{axis}
\end{tikzpicture}
	\caption{Solution times of the linear relaxation for 14 test cases with $T=120$ periods and $\StartupTol = 5\%$ in the German power system (upper chart) and the IEEE 118 bus system (lower chart). Temp outperforms 3-Bin consistently.}
	\label{figure:IntGapTime}
\end{figure}
\fi

A criterion for the quality of a formulation is the computational effort for solving its (initial) linear relaxation. To stay as close as possible to the practical application, we tried to remain close to the linear relaxation by disabling the integrality-specific algorithms of the solvers, \ie presolve, integrated cuts, and heuristics. The experiments were conducted with an interior point algorithm, which proved to be significantly faster than the dual simplex across all formulations. Applying the latter decreases the difference between 3-Bin and Temp slightly, while 1-Bin and 1-Bin* are considerably slower.

Fig.~\ref{figure:IntGapTime} compares solution times of the linear relaxations taken over 14~time ranges of length~$\squeeze[0.5]{T=120}$ as described in Section~\ref{section:Scenarios} and using a start-up cost approximation tolerance of~$\squeeze[0.5]{\StartupTol=5\%}$. The results show that Temp significantly outperforms 3-Bin. While 1-Bin and 1-Bin* are on average faster than Temp in the German power system, this is reversed in the IEEE 118 bus system where Temp yields the fastest linear relaxation by a considerable margin.

\ifIntGapVert\else
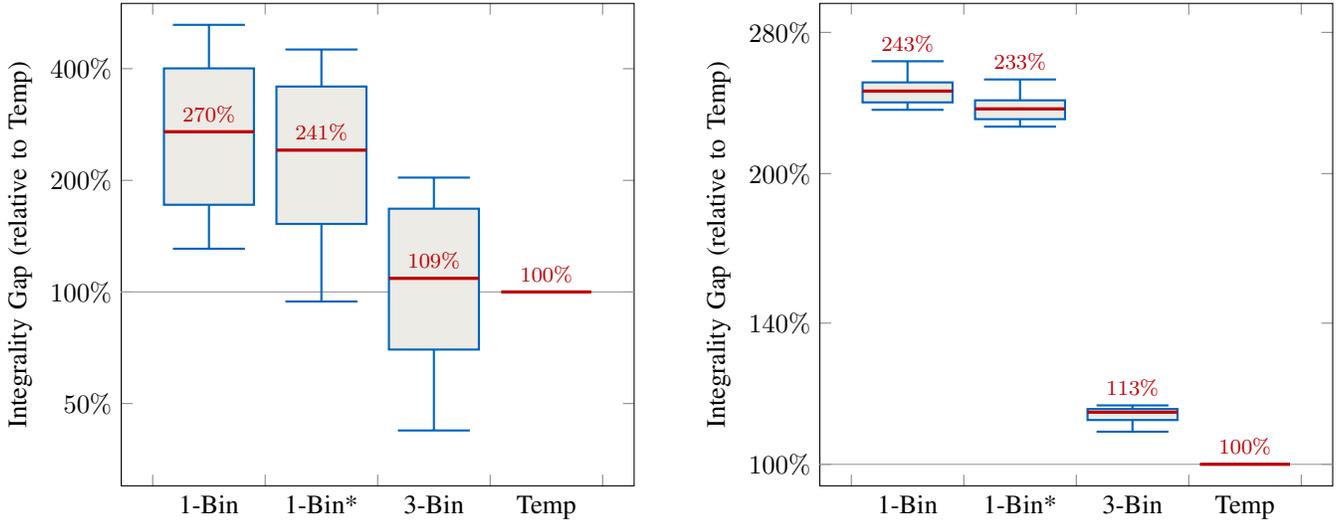
\begin{figure*}[t]
	\hspace*{\fill}
	\begin{tikzpicture}
	\begin{semilogyaxis}[
		myboxplot,
		scaled y ticks = false,
		ylabel={Integrality Gap (relative to Temp)},
		ymin=0.3,
		ymax=6,
		ytick={0.5,1,2,4},
		yticklabels={$50\%$,$100\%$,$200\%$,$400\%$}
	]

	\draw[black!40] (axis cs:0,1) -- (axis cs:6,1);

\addplot[mark=+,boxplot prepared={
    lower whisker=1.307868371475805,
    lower quartile=1.718466764848932,
    median=2.704145556587299,
    upper quartile=4.008535919443189,
    upper whisker=5.252550921109848
}] table[row sep=\\,y index=0] {data\\};
\node[anchor=south,medianColor] at (axis cs:1,2.704145556587299) {\footnotesize$270\%$};
% printBoxplot([1.313608359758285 2.005663885623598 4.906718725303688 1.690767625087371 1.718466764848932 4.276110748202392 2.726878285686959 2.681412827487639 3.466193797008839 4.008535919443189 5.252550921109848 3.830591210213919 2.501542979211559 1.307868371475805])
\addplot[mark=+,boxplot prepared={
    lower whisker=0.9423994540557632,
    lower quartile=1.525747045693449,
    median=2.412329724265502,
    upper quartile=3.581364381201154,
    upper whisker=4.506554439573497
}] table[row sep=\\,y index=0] {data\\};
\node[anchor=south,medianColor] at (axis cs:2,2.412329724265502) {\footnotesize$241\%$};
% printBoxplot([0.9423994540557632  1.681914576635593  4.251538715567866  1.425616531025003  1.525747045693449  3.917703166837309  2.457642218048231  2.367017230482773  2.937957965254645  3.581364381201154  4.506554439573497  3.455373935179408  2.187816003360271  1.164884180623059])
\addplot[mark=+,boxplot prepared={
    lower whisker=0.4231026182180462,
    lower quartile=0.6993244189937931,
    median=1.088209932085982,
    upper quartile=1.677148504035497,
    upper whisker=2.034335022068329
}] table[row sep=\\,y index=0] {data\\};
\node[anchor=south,medianColor] at (axis cs:3,1.088209932085982) {\footnotesize$109\%$};
% printBoxplot([0.4231026182180462 0.7336651719131625  1.677148504035497 0.8225499814419431 0.5068179397226357  1.900966983412628  1.277288198431575  1.031420497888917  1.144999366283047  1.543378980635627  2.034335022068329  1.839109000890271 0.6993244189937931 0.6005780059130617])
\addplot[mark=+,boxplot prepared={
    lower whisker=1,
    lower quartile=1,
    median=1,
    upper quartile=1,
    upper whisker=1
}] table[row sep=\\,y index=0] {data\\};
\node[anchor=south,medianColor] at (axis cs:4,1) {\footnotesize$100\%$};
% printBoxplot([1 1 1 1 1 1 1 1 1 1 1 1 1 1])
    \end{semilogyaxis}
\end{tikzpicture}
	\hspace*{3ex}
	\begin{tikzpicture}
	\begin{semilogyaxis}[
		myboxplot,
		scaled y ticks = false,
		ylabel={Integrality Gap (relative to Temp)},
		ymin=0.95,
		ymax=3,
		ytick={0.7,1,1.4,2,2.8},
		yticklabels={$70\%$,$100\%$,$140\%$,$200\%$,$280\%$}
	]

	\draw[black!40] (axis cs:0,1) -- (axis cs:6,1);

\addplot[mark=+,boxplot prepared={
    lower whisker=2.327925731443056,
    lower quartile=2.369561311316572,
    median=2.434288746063224,
    upper quartile=2.48512930751461,
    upper whisker=2.613638279306984
}] table[row sep=\\,y index=0] {data\\};
\node[anchor=south,medianColor] at (axis cs:1,2.613638279306984) {\footnotesize$243\%$};
% printBoxplot([2.327925731443056 2.369561311316572 2.413924823370401 2.466431276979954  2.48512930751461 2.443648610942137 2.493005009210141  2.43076785022847 2.368402350825332 2.433869052537321 2.362908370692339 2.434708439589128 2.613638279306984 2.565452892856054])
\addplot[mark=+,boxplot prepared={
    lower whisker=2.236355180839924,
    lower quartile=2.276730178371811,
    median=2.333121018338668,
    upper quartile=2.38129694019472,
    upper whisker=2.502606813166213
}] table[row sep=\\,y index=0] {data\\};
\node[anchor=south,medianColor] at (axis cs:2,2.502606813166213) {\footnotesize$233\%$};
% printBoxplot([2.236355180839924 2.276730178371811 2.310380966118722 2.360953299288934  2.38129694019472 2.343907502507977  2.38296501184057 2.325413481184422 2.268460672413853 2.331227088592927 2.262760247714257 2.335014948084408 2.502606813166213 2.448019801665327])
\addplot[mark=+,boxplot prepared={
    lower whisker=1.080846758648522,
    lower quartile=1.111363224080058,
    median=1.132320273999312,
    upper quartile=1.140979633078147,
    upper whisker=1.150784534736868
}] table[row sep=\\,y index=0] {data\\};
\node[anchor=south,medianColor] at (axis cs:3,1.150784534736868) {\footnotesize$113\%$};
% printBoxplot([1.102146297703335 1.130411384722746 1.130972434464579 1.140979633078147  1.14206918884717 1.111363224080058  1.14017843358215 1.134307556982975  1.11881521768799 1.133668113534044 1.080846758648522 1.101582637278001 1.146293800471516 1.150784534736868])
\addplot[mark=+,boxplot prepared={
    lower whisker=1,
    lower quartile=1,
    median=1,
    upper quartile=1,
    upper whisker=1
}] table[row sep=\\,y index=0] {data\\};
\node[anchor=south,medianColor] at (axis cs:4,1) {\footnotesize$100\%$};
% printBoxplot([1 1 1 1 1 1 1 1 1 1 1 1 1 1])
    \end{semilogyaxis}
\end{tikzpicture}
	\hspace*{\fill}
	\caption{Integrality gaps relative to Temp, for 14 test cases with~$T=72$ periods in the German power system (left chart) and the IEEE 118 bus system (right chart). In both, 3-Bin dominates 1-Bin and 1-Bin*, but is on average inferior to Temp. Results in the IEEE 118 bus system exhibit less variance.}
	\label{figure:IntGapRel}
\end{figure*}
\fi

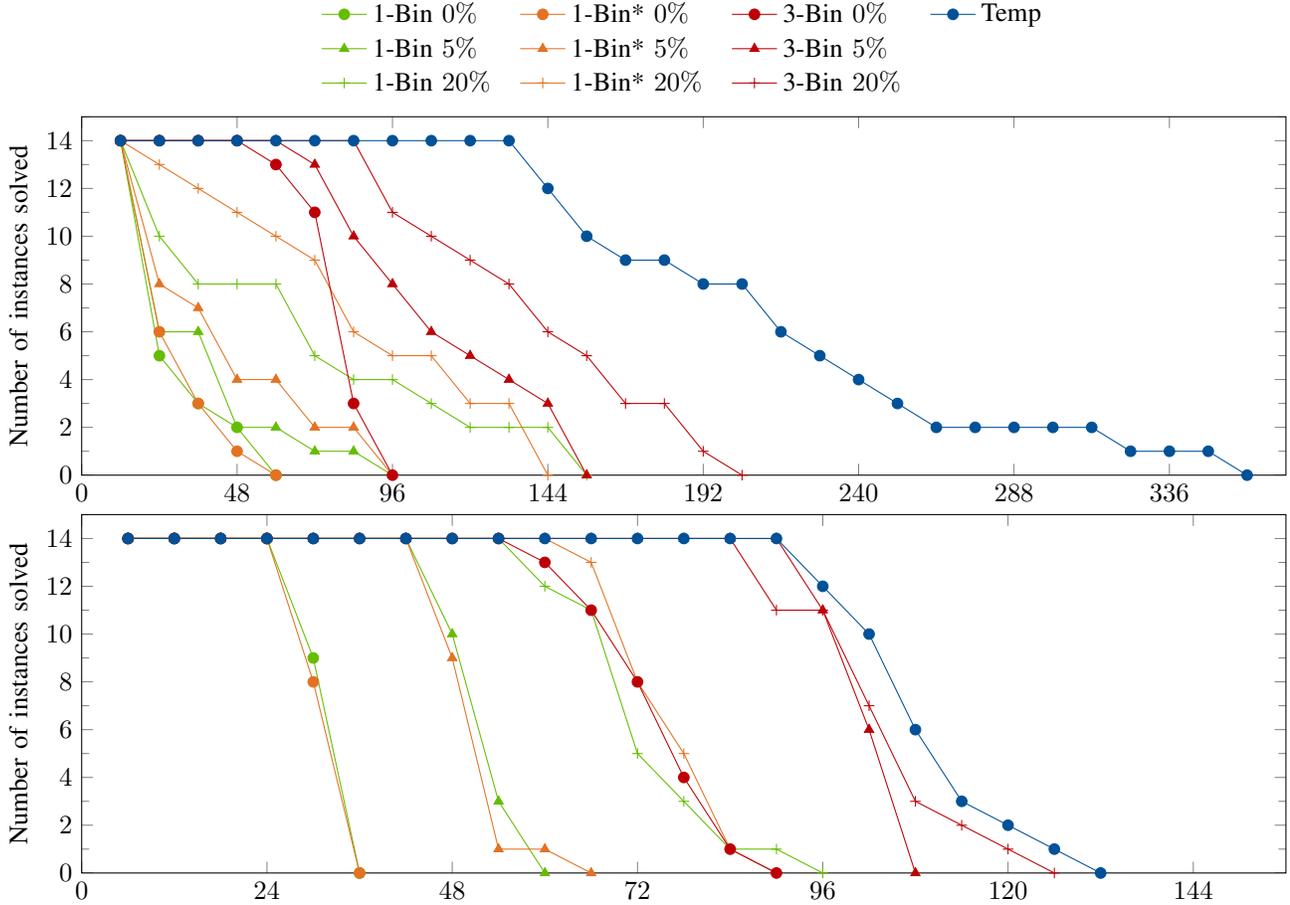
\begin{figure*}[t]
	\centering
	   \pgfplotsset{
	      yticklabel style={overlay},
	      ylabel style={overlay},
	   }% 
	\begin{tikzpicture}
	\begin{axis}[
		width=0.97\textwidth,
		height=6.35cm,
		xmin=0,
		xmax=372,
		xtick={0,48,96,144,192,240,288,336},
		ymin=0,
		ymax=15,
		restrict y to domain=0:14,
		ytick={0,2,4,6,8,10,12,14},
		minor y tick num=1,
		ylabel={Number of instances solved},
		legend columns=3,
		transpose legend,
		legend style={at={(0.5,1.03)},anchor=south,draw=none,fill=white,legend cell align=left,/tikz/every even column/.append style={column sep=2ex}},
	]
	
		\def\file{pics/nrsolved_DE.csv}
	
	\addplot[mark=*,		color=1Bin]	 table[x=nrSteps, y=1Bin0, col sep=comma]  {\file};
		\addplot[mark=triangle*,color=1Bin]	 table[x=nrSteps, y=1Bin5, col sep=comma]  {\file};
		\addplot[mark=+,		color=1Bin]  table[x=nrSteps, y=1Bin20, col sep=comma]  {\file};
				
		\addplot[mark=*,		color=1BinS] table[x=nrSteps, y=1Bins0, col sep=comma]  {\file};
		\addplot[mark=triangle*,color=1BinS] table[x=nrSteps, y=1Bins5, col sep=comma]  {\file};
		\addplot[mark=+,		color=1BinS] table[x=nrSteps, y=1Bins20, col sep=comma]  {\file};
		
		\addplot[mark=*,		color=3Bin]  table[x=nrSteps, y=3Bin0, col sep=comma]  {\file};
		\addplot[mark=triangle*,color=3Bin]  table[x=nrSteps, y=3Bin5, col sep=comma] {\file};
		\addplot[mark=+,		color=3Bin]  table[x=nrSteps, y=3Bin20, col sep=comma] {\file};
	
		\addplot[mark=*,		color=Tbase]  table[x=nrSteps, y=Temp, col sep=comma]  {\file};
		%\addplot[mark=diamond*,	color=Temp] table[x=nrSteps, y=T-Cuts, col sep=comma]  {\file};
				\addlegendimage{empty legend}		
		
				\addlegendentry{1-Bin~$0\%$};	
				\addlegendentry{1-Bin~$5\%$};
				\addlegendentry{1-Bin~$20\%$};
		
				\addlegendentry{1-Bin*~$0\%$};
				\addlegendentry{1-Bin*~$5\%$};
				\addlegendentry{1-Bin*~$20\%$};
		
				\addlegendentry{3-Bin~$0\%$};
				\addlegendentry{3-Bin~$5\%$};
				\addlegendentry{3-Bin~$20\%$};
		
				\addlegendentry{Temp};
				\addlegendentry{};

	\end{axis}
\end{tikzpicture}
	\vskip0.5mm
	\begin{tikzpicture}
	\begin{axis}[
width=0.97\textwidth,
		height=6.35cm,
		xmin=0,
		xmax=156,
		xtick={0,24,48,72,96,120,144},
		ymin=0,
		ymax=15,
		restrict y to domain=0:14,
		ytick={0,2,4,6,8,10,12,14},
		minor y tick num=1,
		ylabel={Number of instances solved},
	]
	
		\def\file{pics/nrsolved_IEEE118.csv}
	
	\addplot[mark=*,		color=1Bin]	 table[x=nrSteps, y=1Bin0, col sep=comma]  {\file};
		\addplot[mark=triangle*,color=1Bin]	 table[x=nrSteps, y=1Bin5, col sep=comma]  {\file};
		\addplot[mark=+,		color=1Bin]  table[x=nrSteps, y=1Bin20, col sep=comma]  {\file};
				
		\addplot[mark=*,		color=1BinS] table[x=nrSteps, y=1Bins0, col sep=comma]  {\file};
		\addplot[mark=triangle*,color=1BinS] table[x=nrSteps, y=1Bins5, col sep=comma]  {\file};
		\addplot[mark=+,		color=1BinS] table[x=nrSteps, y=1Bins20, col sep=comma]  {\file};
		
		\addplot[mark=*,		color=3Bin]  table[x=nrSteps, y=3Bin0, col sep=comma]  {\file};
		\addplot[mark=triangle*,color=3Bin]  table[x=nrSteps, y=3Bin5, col sep=comma] {\file};
		\addplot[mark=+,		color=3Bin]  table[x=nrSteps, y=3Bin20, col sep=comma] {\file};
	
		\addplot[mark=*,		color=Tbase]  table[x=nrSteps, y=Temp, col sep=comma]  {\file};
		%\addplot[mark=diamond*,	color=Temp] table[x=nrSteps, y=T-Cuts, col sep=comma]  {\file};

	\end{axis}
\end{tikzpicture}
	\caption[Scaling of computational effort with problem size]{Scaling of computational effort with problem size for all formulations and different start-up cost approximation tolerances~$\StartupTol$. The upper chart shows the number of instances solved to an optimality gap of 1\% within 30~minutes of computation time for the German power system, the lower chart shows the same for the IEEE 118 bus system.}
	\label{figure:scalingcomparesolved}
\end{figure*}

\subsection{Integrality Gap}
\label{section:IntegralityGap}

Another important criterion of a problem formulation is its integrality gap, which measures the influence of the integrality constraints on the optimal solution. It is defined as $z_{\textup{MIP}}-z_{\textup{LP}}$, where $z_{\textup{MIP}}$ denotes the optimal value and $z_{\textup{LP}}$ the optimal fractional value, and normalized to $\sfrac{(z_{\textup{MIP}}-z_{\textup{LP}})}{z_{\textup{MIP}}}$ for comparability across test cases.

Smaller integrality gaps mean better lower bounds, which lead to faster solution times. The best possible integrality gap is~$0$, which would mean that the optimal objective value of the formulation does not depend on the integrality constraints.

Fig.~\ref{figure:IntGapRel} shows the integrality gap of all four models relative to Temp for the same test cases as in Fig.~\ref{figure:IntGapTime}, but with $\StartupTol = 0\%$ and $T = 72$. Note that the medians of the relative integrality gaps are similar in both formulations, but the German power system exhibits a much higher variance. We attribute this to the highly volatile residual demand in our forecast of the year 2025. Moreover, the complexity of the extended formulation for the IEEE 118 bus system results in a higher absolute integrality gap, with a median of $1.4\%$ compared to $0.7\%$ in the basic formulation.

Fig.~\ref{figure:IntGapRel} clearly illustrates the advantage of modeling the temperature as an explicit variable. Since the inequalities~\eqref{constraint:startcosttight} of 1-Bin* dominate the inequalities~\eqref{constraint:startcostorig} of 1-Bin, 1-Bin* must have a lower integrality gap (11\% and 4\% decrease). 3-Bin consistently provides a lower integrality gap than 1-Bin*, with an average reduction of 55\% and 51\%, corresponding in magnitude to the results in \cite{Morales-Espana_Thermal_2013}. Temp further decreases the average integrality gap of 3-Bin by 9 and 13 percentage points and proves to be the tightest model analyzed. 

\ifIntGapVert
\begin{figure}[b]
	\begin{tikzpicture}
	\begin{semilogyaxis}[
		myboxplot,
		scaled y ticks = false,
		ylabel={Integrality Gap (relative to Temp)},
		ymin=0.3,
		ymax=6,
		ytick={0.5,1,2,4},
		yticklabels={$50\%$,$100\%$,$200\%$,$400\%$}
	]

	\draw[black!40] (axis cs:0,1) -- (axis cs:6,1);

\addplot[mark=+,boxplot prepared={
    lower whisker=1.307868371475805,
    lower quartile=1.718466764848932,
    median=2.704145556587299,
    upper quartile=4.008535919443189,
    upper whisker=5.252550921109848
}] table[row sep=\\,y index=0] {data\\};
\node[anchor=south,medianColor] at (axis cs:1,2.704145556587299) {\footnotesize$270\%$};
% printBoxplot([1.313608359758285 2.005663885623598 4.906718725303688 1.690767625087371 1.718466764848932 4.276110748202392 2.726878285686959 2.681412827487639 3.466193797008839 4.008535919443189 5.252550921109848 3.830591210213919 2.501542979211559 1.307868371475805])
\addplot[mark=+,boxplot prepared={
    lower whisker=0.9423994540557632,
    lower quartile=1.525747045693449,
    median=2.412329724265502,
    upper quartile=3.581364381201154,
    upper whisker=4.506554439573497
}] table[row sep=\\,y index=0] {data\\};
\node[anchor=south,medianColor] at (axis cs:2,2.412329724265502) {\footnotesize$241\%$};
% printBoxplot([0.9423994540557632  1.681914576635593  4.251538715567866  1.425616531025003  1.525747045693449  3.917703166837309  2.457642218048231  2.367017230482773  2.937957965254645  3.581364381201154  4.506554439573497  3.455373935179408  2.187816003360271  1.164884180623059])
\addplot[mark=+,boxplot prepared={
    lower whisker=0.4231026182180462,
    lower quartile=0.6993244189937931,
    median=1.088209932085982,
    upper quartile=1.677148504035497,
    upper whisker=2.034335022068329
}] table[row sep=\\,y index=0] {data\\};
\node[anchor=south,medianColor] at (axis cs:3,1.088209932085982) {\footnotesize$109\%$};
% printBoxplot([0.4231026182180462 0.7336651719131625  1.677148504035497 0.8225499814419431 0.5068179397226357  1.900966983412628  1.277288198431575  1.031420497888917  1.144999366283047  1.543378980635627  2.034335022068329  1.839109000890271 0.6993244189937931 0.6005780059130617])
\addplot[mark=+,boxplot prepared={
    lower whisker=1,
    lower quartile=1,
    median=1,
    upper quartile=1,
    upper whisker=1
}] table[row sep=\\,y index=0] {data\\};
\node[anchor=south,medianColor] at (axis cs:4,1) {\footnotesize$100\%$};
% printBoxplot([1 1 1 1 1 1 1 1 1 1 1 1 1 1])
    \end{semilogyaxis}
\end{tikzpicture}
	\vskip1ex
	\begin{tikzpicture}
	\begin{semilogyaxis}[
		myboxplot,
		scaled y ticks = false,
		ylabel={Integrality Gap (relative to Temp)},
		ymin=0.95,
		ymax=3,
		ytick={0.7,1,1.4,2,2.8},
		yticklabels={$70\%$,$100\%$,$140\%$,$200\%$,$280\%$}
	]

	\draw[black!40] (axis cs:0,1) -- (axis cs:6,1);

\addplot[mark=+,boxplot prepared={
    lower whisker=2.327925731443056,
    lower quartile=2.369561311316572,
    median=2.434288746063224,
    upper quartile=2.48512930751461,
    upper whisker=2.613638279306984
}] table[row sep=\\,y index=0] {data\\};
\node[anchor=south,medianColor] at (axis cs:1,2.613638279306984) {\footnotesize$243\%$};
% printBoxplot([2.327925731443056 2.369561311316572 2.413924823370401 2.466431276979954  2.48512930751461 2.443648610942137 2.493005009210141  2.43076785022847 2.368402350825332 2.433869052537321 2.362908370692339 2.434708439589128 2.613638279306984 2.565452892856054])
\addplot[mark=+,boxplot prepared={
    lower whisker=2.236355180839924,
    lower quartile=2.276730178371811,
    median=2.333121018338668,
    upper quartile=2.38129694019472,
    upper whisker=2.502606813166213
}] table[row sep=\\,y index=0] {data\\};
\node[anchor=south,medianColor] at (axis cs:2,2.502606813166213) {\footnotesize$233\%$};
% printBoxplot([2.236355180839924 2.276730178371811 2.310380966118722 2.360953299288934  2.38129694019472 2.343907502507977  2.38296501184057 2.325413481184422 2.268460672413853 2.331227088592927 2.262760247714257 2.335014948084408 2.502606813166213 2.448019801665327])
\addplot[mark=+,boxplot prepared={
    lower whisker=1.080846758648522,
    lower quartile=1.111363224080058,
    median=1.132320273999312,
    upper quartile=1.140979633078147,
    upper whisker=1.150784534736868
}] table[row sep=\\,y index=0] {data\\};
\node[anchor=south,medianColor] at (axis cs:3,1.150784534736868) {\footnotesize$113\%$};
% printBoxplot([1.102146297703335 1.130411384722746 1.130972434464579 1.140979633078147  1.14206918884717 1.111363224080058  1.14017843358215 1.134307556982975  1.11881521768799 1.133668113534044 1.080846758648522 1.101582637278001 1.146293800471516 1.150784534736868])
\addplot[mark=+,boxplot prepared={
    lower whisker=1,
    lower quartile=1,
    median=1,
    upper quartile=1,
    upper whisker=1
}] table[row sep=\\,y index=0] {data\\};
\node[anchor=south,medianColor] at (axis cs:4,1) {\footnotesize$100\%$};
% printBoxplot([1 1 1 1 1 1 1 1 1 1 1 1 1 1])
    \end{semilogyaxis}
\end{tikzpicture}
	\caption{Integrality gaps relative to 3-Bin, for 14 test cases with~$T=72$ periods in the German power system (upper chart) and the IEEE 118 bus system (lower chart). In both, 3-Bin dominates 1-Bin and 1-Bin*, but is on average inferior to Temp. Results in the IEEE 118 bus system exhibit less variance.}
	\label{figure:IntGapRel}
\end{figure}
\fi

\subsection{Performance With Scaling to a Larger Number of Periods}
\label{section:Scaling}

An essential aspect in computational efficiency is the behavior with model scaling. Especially scaling to a greater number of modeled periods seems to be highly relevant for future operational planning for two reasons:
\begin{enumerate}
  \item As the residual load will become more volatile it will be beneficial to increase the time resolution \cite{Deane2014}.
  \item As renewable generation changes over several days and weeks, the storage management requires to consider longer time horizons than today where it is mainly driven by day and night variation of load.
\end{enumerate}

We analyze the scaling behavior by considering 14 sets of test cases with start periods~$S$ described in Section~\ref{section:Scenarios}, $T$ varying from $\squeeze[0.5]{T=24}$ to $\squeeze[0.5]{T=444}$, and start-up costs approximated to tolerances $\squeeze[0.8]{\StartupTol \in \{0\%,5\%,20\%\}}$. Fig.~\ref{figure:scalingcomparesolved} shows the number of instances which have been solved to an optimality gap of $1\%$ within 30~minutes for the German power sytem (upper chart) and the IEEE 118 bus system (lower chart).

In both cases, 3-Bin dominates 1-Bin and 1-Bin*, confirming the results in \cite{Morales-Espana_Thermal_2013}. However, even if we allow the highest start-up cost approximation tolerance $\StartupTol = 20\%$, Temp consistently solves more instances than all other models.

Unsurprisingly, the superiority of the temperature model is more emphasized in the basic formulation, since the higher complexity of the extended formulation and the lower share of start-up costs in the IEEE 118 system lessen the impact of the start-up cost model.  % Numerical examples
\section{Conclusion}
\label{section:Conclusion}

Firstly, we presented the tightened formulation 1-Bin* of one of the state-of-the-art models for general increasing start-up cost functions. The main result however is the temperature model, which accurately models the exponential start-up costs while consistently outperforming the existing formulations---even when allowing high approximation tolerances. The increased performance should be attributed mainly to the typically smaller integrality gap.

We hope that the novel temperature formulation and its physical interpretation will inspire and facilitate fundamental future extensions of the Unit Commitment problem.  % Conclusion
\appendices
\section{Ramping and Minimum Up-/Down Constraints}
\label{section:Appendix}

\noindent The ramping speed of a unit is described by four parameters,
\begin{itemize}
  \item the maximum ramp up speed~$\rampup(i)$ when operational,
  \item the maximum ramp down speed~$\rampdown(i)$ when operational,
  \item the maximum ramp up at start-up~$\startupRamp(i)$, and
  \item the maximum ramp down at shutdown~$\shutdownRamp(i)$.
\end{itemize}

\noindent These ramping speeds are modeled in \cite{Carrion2006} as
\begin{IEEEeqnarray}{rCl}
	\label{constraint:RampUp}
	\prod(i,t) &\leq& \prod(i,t-1) + \rampup(i) \onOff(i,t-1) + \startupRamp(i) (\onOff(i,t) - \onOff(i,t-1)) + \maxProd(i)(1-\onOff(i,t))\nonumber\\
		\IEEEeqnarraymulticol{3}{r}{\forall\,i \in \I, t \in \discrange{2}{T},\IEEEeqnarraynumspace}\\[1ex]
	\label{constraint:RampDown}
	\prod(i,t) &\geq& \prod(i,t-1) - \rampdown(i) \onOff(i,t) - \shutdownRamp(i) (\onOff(i,t-1) - \onOff(i,t)) - \maxProd(i)(1-\onOff(i,t-1)) \IEEEnonumber\\
		\IEEEeqnarraymulticol{3}{r}{\forall\,i \in \I, t \in \discrange{2}{T},\IEEEeqnarraynumspace}\\[1ex]
	\label{constraint:Shutdown}
	\prod(i,t) &\leq& \maxProd(i) \onOff(i,t+1) + \shutdownRamp(i) (\onOff(i,t) - \onOff(i,t+1)) \IEEEnonumber\\
		\IEEEeqnarraymulticol{3}{r}{\forall\,i \in \I,t \in \discrange*{1}{T-1}.\IEEEeqnarraynumspace}
\end{IEEEeqnarray}

In \cite{Ostrowski2012}, a tighter version of the ramping constraints is proposed, which uses start-up and shutdown indicators. Its main inequalities are
\begin{IEEEeqnarray}{rCll}
	\label{equation:OstrowskiStart}
	% Inequality (6)
	\prod(i,t) - \prod(i,t-1)
	&\leq& \rampup(i) \onOff(i,t-1)
	+ \startupRamp(i)\startupIndicator(i,t)
	\qquad &\forall\,i \in \I, t \in \discrange{2}{T},\IEEEeqnarraynumspace\\[1ex]
%
	% Inequality (7)
	\label{equation:Ostrowski2}
	\prod(i,t-1) - \prod(i,t)
	&\leq& \rampdown(i) \onOff(i,t)
	+ \shutdownRamp(i) \shutdownIndicator(i,t)
	\qquad &\forall\,i \in \I, t \in \discrange{2}{T}.\IEEEeqnarraynumspace
\end{IEEEeqnarray}

\noindent Depending the unit parameters, \cite{Ostrowski2012} furthermore adds:
\begin{itemize}
  \item if $\rampdown(i)>(\startupRamp(i)-\minProd(i))$ and $\minup(i) \geq 2$,
	\begin{IEEEeqnarray}{rCl}
		% Inequality (20)
		\label{equation:Ostrowski3}
		\prod(i,t-1) - \prod(i,t)
		&\leq& \rampdown(i) \onOff(i,t)
		+ \shutdownRamp(i) \shutdownIndicator(i,t)  \nonumber\\
		&&- (\rampdown(i) - \startupRamp(i) + \minProd(i)) \startupIndicator(i,t-1)\\
		&&-(\rampdown(i) + \minProd(i)) \startupIndicator(i,t) \qquad \forall\, t \in \discrange{2}{T}, \nonumber
	\end{IEEEeqnarray}
	\item if $\rampdown(i) > (\startupRamp(i)-\minProd(i)), \minup(i) \geq 3$ and $\mindown(i) \geq 2$,
	\begin{IEEEeqnarray}{rCl}
	% Inequality (21)
		\label{equation:Ostrowski4}
		\prod(i,t-1)-\prod(i,t)
		&\leq& \rampdown(i) \onOff(i,t+1)
		- (\rampdown(i) - \startupRamp(i) + \minProd(i)) \startupIndicator(i,t-1)\IEEEeqnarraynumspace\nonumber\\
		&&- (\rampdown(i) + \minProd(i)) \startupIndicator(i,t)
		- \rampdown(i) \startupIndicator(i,t+1)\\
		&&+ \shutdownRamp(i) \shutdownIndicator(i,t)
		+ \rampdown(i) \shutdownIndicator(i,t+1)
		\forall\,t \in \discrange*{2}{T-1},\nonumber
	\end{IEEEeqnarray}
	\item for all units,
	\begin{IEEEeqnarray}{rCl}
		% Inequality (22)
		\label{equation:Ostrowski5}
		\prod(i,t-2) - \prod(i,t)
		&\leq& 2\rampdown(i) \onOff(i,t)
		+ \shutdownRamp(i) \shutdownIndicator(i,t-1)
		+ (\shutdownRamp(i) + \rampdown(i))\shutdownIndicator(i,t)\nonumber\\
		&&- 2 \rampdown(i)\startupIndicator(i,t-2)
		- (2\rampdown(i) + \minProd(i))	\startupIndicator(i,t-1)\\
		&&- (2\rampdown(i) + \minProd(i))\startupIndicator(i,t) \quad \forall\,i \in \I, t \in \discrange*{3}{T},\nonumber
	\end{IEEEeqnarray}
	\item if $\rampup(i)>(\shutdownRamp(i)-\minProd(i))$  and $\minup(i) \geq 2$,
	\begin{IEEEeqnarray}{rCl}
		% Inequality (23)
		\label{equation:Ostrowski6}
		\prod(i,t) - \prod(i,t-1)
		&\leq& \rampup(i) \onOff(i,t) + (\startupRamp(i) - \rampup(i)) \startupIndicator(i,t)\nonumber\\
		&&- \minProd(i)\shutdownIndicator(i,t) - (\rampup(i) - \shutdownRamp(i) + \minProd(i)) \shutdownIndicator(i,t+1)\hspace*{2em}\\
		&&\hspace{4em}\qquad \forall\,i \in \I, t \in \discrange*{2}{T-1},\nonumber
	\end{IEEEeqnarray}
	\item if $\rampup(i)>(\shutdownRamp(i)-\minProd(i))$ and $\mindown(i) \geq 2$,
	\begin{IEEEeqnarray}{rCl}
		% Inequality (24)
		\label{equation:OstrowskiEnd}
		\prod(i,t) - \prod(i,t-2)
		&\leq& 2 \rampup(i) \onOff(i,t)
		- \minProd(i) \shutdownIndicator(i,t-1)
		- \minProd(i) \shutdownIndicator(i,t)\nonumber\\
		&&+(\startupRamp(i) - \rampup(i)) \startupIndicator(i,t-1)
		+ (\startupRamp(i) - 2\rampup(i)) \startupIndicator(i,t)\hspace*{2em}\\
		&&\qquad \forall\,i \in \I, t \in \discrange*{3}{T-1}.\nonumber
	\end{IEEEeqnarray}
\end{itemize}

The minimum up-/downtime is modeled in \cite{Rajan2005} by the turn on/off inequalities,
\begin{alignat}{2}
\label{equation:minup1}
	\sum_{k=t-\minup(i)+1}^{t} \startupIndicator(i,k) &\leq \onOff(i,t) &&\qquad \forall\,i \in \I, t \in \discrange{\minup(i)}{T},\\
\label{equation:minup2}
	\sum_{k=t-\mindown(i)+1}^{t} \shutdownIndicator(i,k) &\leq 1-\onOff(i,t) &&\qquad \forall\,i \in \I, t \in \discrange{\mindown(i)}{T}.
\end{alignat}

% Bibliography
\ifCLASSOPTIONcaptionsoff
  \newpage
\fi

\bibliographystyle{IEEEtran}
\bibliography{refs}

\newpage

% Biographies

\begin{IEEEbiography}[{\includegraphics[width=1in,height=1.25in,clip,keepaspectratio]{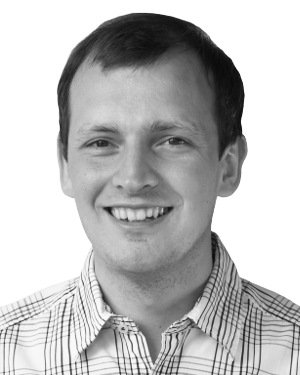}}]{Matthias Silbernagl}
received his Dipl.-Math. (M.Sc. in Mathematics) from the Technische Uni\-ver\-si\-t\"at M\"un\-chen (TUM), Germany, in 2009. Currently, he is a doctoral candidate at the Chair for Applied Geometry and Discrete Mathematics at TUM. As an IGSSE associate, he participates in the interdisciplinary project \enquote{Integration of Renewables}.

His area of research includes Mixed Integer Programming and Polyhedral Studies, with a focus on the family of Unit Commitment problems.
\end{IEEEbiography}

\begin{IEEEbiography}[{\includegraphics[width=1in,height=1.25in,clip,keepaspectratio]{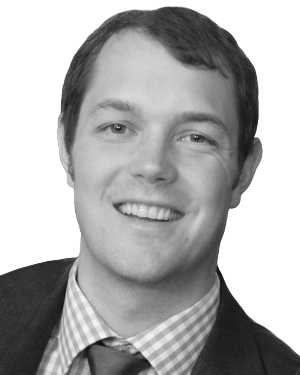}}]{Matthias Huber}
received his Dipl.-Ing. (M.Sc.) in Mechanical Engineering from Technische Universit\"at M\"unchen (TUM) and his B.Sc. in Economics from Ludwig Maximilians Universit\"at (LMU), both in 2010. He is now pursuing his PhD at the Institute for Renewable and Sustainable Energy Systems at TUM. His research interests include optimal planning, operation, and economics of power systems.
\end{IEEEbiography}

\begin{IEEEbiography}[{\includegraphics[width=1in,height=1.25in,clip,keepaspectratio]{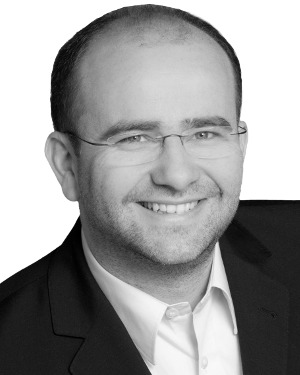}}]{Ren\'e Brandenberg}
received his Dipl-Math. from Uni\-ver\-si\-t\"at Trier in 1998 and his Dr. rer. nat. (Ph.D.) from TUM in 2003. He is now on a permanent position
at the Chair for Applied Geometry and Discrete Mathematics at TUM and his main research areas are convex and computational geometry as well as their application in optimization.
\end{IEEEbiography}

\vfill

\end{document}